\documentclass[11pt]{article}
\usepackage[english]{babel}
\usepackage[cp850]{inputenc}
\usepackage[dvips]{graphicx}
\usepackage{amsmath,amsfonts,amsthm,amssymb}
\usepackage[usenames, dvipsnames]{color}
\usepackage{fancyhdr}
\usepackage{stmaryrd}
\usepackage[colorlinks=true,citecolor=red,linkcolor=blue,urlcolor=blue,pdfstartview=FitH]{hyperref}
\usepackage{dsfont}
\usepackage{xcolor}
\usepackage{epsfig}

\font\tenmath=msbm10 scaled 1200

\font\sevenmath=msbm7 scaled 1200
\font\fivemath=msbm5 scaled 1200 

\newfam\mathfam \textfont\mathfam=\tenmath
\scriptfont\mathfam=\sevenmath \scriptscriptfont\mathfam=\fivemath

\def\R{{\mathbb{R}}}
\def\N{{\mathbb{N}}}
\def\E{{\mathbb{E}}}
\def\L{{\cal L}}
\def\A{{\cal A}}
\def\P{{\mathbb{P}}}

\def\C{{\mathbb{C}}}
\def\F{{\cal F}}

\newtheorem{theo}{Theorem}[section]
\newtheorem{lem}{Lemma}[section]

\newtheorem{cor}{Corollary}[section]

\newtheorem{pro}{Proof}[section]

\newfam\mathfam \textfont\mathfam=\tenmath
\scriptfont\mathfam=\sevenmath \scriptscriptfont\mathfam=\fivemath

\def \^#1{\if#1i{\accent"5E\i}\else{\accent"5E#1}\fi}

\def \cadre{(\Omega, {\cal A}, \P)}

\def \cqfd{\quad\Box}

\def \bs{\bigskip}
\def \ni{\noindent}
\def \Tr{w}

\title{\bf Randomized Urn Models revisited using Stochastic Approximation}
 
\author{\textsc{Sophie Laruelle} \thanks{Laboratoire de Probabilit\'es et Mod\`eles al\'eatoires, UMR~7599, Universit\'e Paris 6, case 188, 4, pl. Jussieu, F-75252 Paris Cedex 5, France. E-mail: \texttt{sophie.laruelle@upmc.fr}}  
\and \textsc{Gilles Pag\`es} \thanks{Laboratoire de Probabilit\'es et Mod\`eles al\'eatoires, UMR~7599, Universit\'e Paris 6
. E-mail: \texttt{gilles.pages@upmc.fr}}}

\date{September 13, 2016}

\begin{document}

\maketitle 

\vspace{-1cm}


\vspace{0.5cm}

\begin{abstract}
This paper presents the link between stochastic approximation and clinical trials based on randomized urn models investigated in~\cite{BaiHu,BaiHu2,BaiHuShe}. We reformulate the  dynamics of both the urn composition and the  assigned  treatments  as  standard stochastic approximation ($SA$) algorithms with remainder. Then, we derive the $a.s.$ convergence and the  asymptotic normality (Central Limit Theorem $CLT$) of the normalized procedure  under less stringent assumptions by calling upon the  $ODE$  and $SDE$ methods. As a second step, we investigate a more involved family of models, known as multi-arm clinical trials, where the urn updating depends on the past performances of the treatments. By increasing the dimension of the state vector, our $SA$ approach provides this time a new asymptotic normality  result.
\end{abstract}

\bigskip

\small
\noindent {\em This is the extended version of the eponym published paper in {\em Annals of Applied Probability} {\bf 23}(4):1409-1436. Proofs are more detailed and additional results are established on specified models of urns investigated in the paper. }


\bs \ni {\em 2010 AMS classification:} 62L20, 62E20, 62L05
secondary: 62F12, 62P10.

\section{Introduction}

The aim of this paper is to illustrate the efficiency of Stochastic Approximation ($SA$) Theory by revisiting several recent results on randomized urn models applied to clinical trials (especially~\cite{BaiHu,BaiHu2,BaiHuShe}). We will first retrieve the $a.s.$ convergence (strong consistency) and asymptotic normality results obtained in these papers under less stringent assumptions. Then we will take advantage of this more synthetic approach to establish a new Central Limit Theorem ($CLT$) in the more sophisticate randomized urn model known as ``multi-arm clinical test". In this model, the urn updating which produces the adaptive design is based on statistical estimators  of the past efficiency of the assigned treatments. 

In these adaptive models, the starting point is the equation which governs  the  urn composition updated after  each new treated patient.  Basically, we will show that a normalized version of this urn composition can be formulated  as a classical recursive stochastic algorithm with  step $\gamma_n=\frac 1n$ which classical Stochastic Approximation Theory deals with. Doing so we will be in position to establish the $a.s.$ convergence of the procedure by calling upon the so-called Ordinary Differential Equation Method ($ODE$ method) and to derive the asymptotic normality  - a $CLT$, to be precise -  from the standard $CLT$ for stochastic algorithms (sometimes called the Stochastic Differential Equation Method ($SDE$ method), see $e.g.$ \cite{Duf2,BMP}). These two main theoretical results   are recalled in a self-contained form in the Appendix. They can be found in all classical textbooks on $SA$ (\cite{BMP}, \cite{Duf}, \cite{Duf2}, \cite{KusYin}) and go back to \cite{KusCla} and~\cite{Bou}. $SA$ Theory is also used in clinical trials to solve dose-finding problems (see for example \cite{Che} and citations therein).

Clinical trials essentially deal with the asymptotic behaviour of the patient allocation to several treatments during the procedure. Adaptive designs in clinical trials aim at detecting ``on line" which treatment should be assigned to more patients, while keeping randomness enough to preserve the basis of treatments. This adaptive approach relies on the cumulative information provided by the responses to treatments of previous patients in order to adjust treatment allocation to the new patients. To this end, many urn models have been suggested in the literature (see \cite{JohKot}, \cite{Zel}, \cite{Wei}, \cite{FloRos} and \cite{Ros}). The most widespread random adaptive model is the Generalized Friedman Urn ($GFU$) (see \cite{AthKar} and more recently \cite{Jan,Pou}), also called Generalized P\'olya Urn ($GPU$). The idea of this modeling is that the urn contains balls of $d$ different types representative of the treatments. All random variables involved in the model are supposed to be defined on the same probability space $\cadre$. Denote $Y_0=(Y_0^i)_{i=1,\ldots,d}\in\R_+^d\setminus\{0\}$ the initial composition of the urn, where $Y_0^i$ denotes the number of balls of type $i$, $i=1,\ldots,d$ (of course a more realistic though not mandatory assumption would be $Y_0\!\in \N^d\setminus\{0\}$). The allocation of the treatments is sequential and the urn composition at draw $n$ is denoted by $Y_n=(Y_n^i)_{i=1,\ldots,d}$. When the $n^{th}$ patient presents, one draws randomly ($i.e.$ uniformly) a ball from the urn with instant replacement. If the ball is of type $j$, then the treatment $j$ is assigned to the $n^{th}$ patient, $j=1,\ldots,d$, $n\geq1$. The urn composition is updated by taking into account the response of the $n^{th}$ patient to the treatment $j$, or the responses of all patients up to the $n^{th}$ one ($i.e.$ the efficiency of the assigned treatment), namely by adding $D^{ij}_n$ balls of type $i$, $i=1,\ldots,d$. The procedure is iterated as long as patients present. Consequently the larger the number of balls of a given type is, the more efficient the treatment is. The urn composition at stage $n$, modeled by an $\R^d$-valued vector $Y_n$, satisfies the following recursive procedure:
\begin{equation}\label{dynamic}
	Y_n=Y_{n-1}+D_nX_n,\quad n\geq1, \quad Y_0\!\in \R_+^d\setminus\{0\},
\end{equation}
with $D_n=(D^{ij}_n)_{1\leq i,j\leq d}$ is the addition rule matrix and $X_n$ is the result of the $n^{th}$ draw and $X_n:(\Omega,\A,\P)\rightarrow\{e^1,\cdots,e^d\}$ models the selected treatment ($\{e^1,\cdots,e^d\}$ denotes the canonical basis of $\R^d$ and $e^j$ stands for treatment $j$). We assume that there is no extinction $i.e.$ $Y_n\!\in \R_+^d\setminus\{0\}$ $a.s.$ for every $n\ge 1$: so is the case if all the entries $D_n^{ij}$ are $a.s.$ non-negative, but other settings can also be taken under consideration (see Section~\ref{deux}). We model the drawing in the urn by setting 
\begin{equation}\label{ConstructX}
X_n=\sum_{j=1}^d\mbox{\bf 1}_{\left\{\frac{\sum_{\ell=1}^{j-1}Y_{n-1}^{\ell}}{\sum_{\ell=1}^dY_{n-1}^\ell }<U_n\leq\frac{\sum_{\ell=1}^{j}Y_{n-1}^\ell}{\sum_{\ell=1}^dY_{n-1}^\ell }\right\}}e^j, \quad n\ge 1,
\end{equation}
where $(U_n)_{n\geq1}$ is i.i.d. with distribution  $U_1\overset{\L}{\sim}{\cal U}_{[0,1]}$. 

Let $\F_n=\sigma(Y_0,U_k,D_k,1\leq k\leq n)$ be the filtration of the procedure. The {\em generating matrices} are defined as the ${\cal F}_n$-compensator of the additions rule sequence $i.e.$ 
$$
H_n=\left(\E\left[D^{ij}_n\,|\,\F_{n-1}\right]\right)_{1\leq i,j\leq d}, \; n\geq1.
$$
Other fields of application can be considered for such procedures like the adaptive asset allocation by an asset manager or a trader. Indeed this has already been done in \cite{LamPagTar} and successfully implemented with multi-armed bandit procedure. Imagine an asset manager who  can  trade  the same financial instrument (tradable asset) on different trading venues. To optimize the  execution of an inventory  of this asset, she can split her orders across these trading destinations. She starts with the initial allocation vector $Y_0$. At stage $n$, she chooses a trading destinations according to the distribution~(\ref{ConstructX}) of $X_n$, then evaluates its performance during one time step and modifies the urn composition (most likely virtually) and proceeds. Thus the normalized urn composition represents the allocation vector among the venues and the addition rule matrices model the successive re-allocations  depending on the past performances of the different trading destinations. 

One may also consider this type of procedure as a strategy to update the composition of a portfolio or even a whole fund, based on the (recent) past performances of the assets.

The first designs under consideration were the homogeneous $GFU$ models where the addition rules $D_n$ are i.i.d. and the so-called generating matrices $H_n=H=\E D_n$ are identical, non-random, with nonnegative entries and irreducible. Hence by the Perron-Frobenius Theorem $H$ has  a unique and positive maximal eigenvalue and an  eigenvector with positive components (see~\cite{AthKar, AthKar2,Fre,Gou}). But the homogeneity of the generating matrix is often not satisfied in practice and inhomogeneous $GFU$ models have been introduced (see \cite{BaiHu}) in which $H_n$ are not random but converge to a deterministic limit $H$, under the assumption that the total number of balls added at each stage is constant. As a third step, the homogeneous Extended P\'olya Urn ($EPU$) models have been introduced in~\cite{Smy} in which only the mean total number of balls added at each stage is constant. This number is called the {\em balance} of the urn and the urn is said {\em balanced}.

Finally, in~\cite{BaiHu2} the authors proposed a nonhomogeneous $EPU$ model because in applications, the addition rule $D_n$ depends on the past history of previous trials (see \cite{AndFarRam}), so that the general generating matrix $H_n$ is usually random. Thus the entries of $H$ may not be all nonnegative ($e.g.$, when there is no replacement after the draw diagonal terms may become negative), and they assume that the matrix $H$ has a unique maximal eigenvalue $\lambda$ with associated (right) eigenvector $v^*=(v^{*,i})_{i=1,\ldots,d}$ with $\sum_{i=1}^dv^{*,i}=1$. Furthermore the conditional expectation of the total number of balls added at each stage was constant.

The first theoretical investigations on these models focused on the asymptotic properties of the urn composition (consistency and asymptotic normality). However, for practical matter, it is clear that  the asymptotic behaviour of the  vector $N_n:=\sum_{k=1}^n X_k$  which stores the treatment allocation among the first $n$  patients is of high interest, especially its variance structure in order to compare several adaptive designs. Thus, in~\cite{BaiHu2} is proved the strong consistency of both (normalized) quantities $Y_n/n$ and $N_n/n$ (under a summability assumption on the generating matrices).

By considering an appropriate  recursive procedure for the normalized urn composition derived from~(\ref{dynamic})
we prove by the $ODE$ method its $a.s.$ convergence toward $v^*$ under a significantly less stringent assumption, namely the minimal requirement that $H_n\overset{a.s.}{\underset{n\rightarrow+\infty}{\longrightarrow}} H$. The $a.s.$ convergence of the treatment allocation frequency $N_n/n$ toward the same $v^*$ follows from the previous one. 

As concerns asymptotic normality, separate  results on these two quantities  are obtained in~\cite{BaiHu2} under an additional assumption on the rate of convergence of the generating matrices $H_n$  toward $H$. On our side we propose to consider a stochastic approximation procedure with remainder satisfied by the higher dimensional vector $(Y_n/n,N_n/n)$. Then, the standard  $CLT$ for $SA$ procedures with remainder directly provides the expected asymptotic normality result for the whole vector  under an  assumption on the $L^2$-rate of convergence of the generating matrices towards their limit (namely $i.e.$ $|\!|\!|H_n-H|\!|\!|=o(n^{-1/2})$) which is again slightly less stringent than the original one. As a result, we   obtain the asymptotic joint  distribution with an  explicit  global covariance structure matrix.

In the end of~\cite{BaiHu2}, an application to multi-arm clinical trials  randomized urn models is proposed. This adaptive design has already  been introduced in~\cite{BaiHuShe} with first consistency results. This kind of models is clearly the most interesting for practitioners since it takes into account the past results of the assigned treatments in the addition rule matrices, denoted $S_n$ at time $n$ ($S^i_n$ denotes the number of cured patients by treatment $i$ among the $N^i_n$ treated ones). The above strong consistency results apply but none of the asymptotic normality works as stated since the generating matrices $H_n$ do not -- in fact {\em cannot} as we will emphasize -- converge at the requested rate. The reason being that they themselves satisfy a $CLT$. However we van overcome this obstacle by increasing once again the structural dimension of the problem: we show that the triplet $(Y_n/n,N_n/n,S_n/n)$ can be written as a recursive $SA$ algorithm with remainder satisfying $a.s.$ convergence and a $CLT$ (provided the limiting generating matrix is still irreducible, etc).  Thus we illustrate on this example that $SA$ Theory is a powerful tool to investigate this kind of adaptive design problem. The main difficulty is to exhibit the appropriate form for the recursion  by making {\em a priori} the balance between significant asymptotic terms and remainder terms. 

\medskip
The paper is organized as follows. We rewrite the dynamics~(\ref{dynamic}) of the urn composition as a stochastic approximation procedure with state variable for $\widetilde{Y}_n:=Y_n/n$ in Section~\ref{deux1}. In Section~\ref{deux2} the $a.s.$ convergence of $\frac{1}{n}\sum_{i=1}^dY_n^i$ is established which implies that of $\widetilde{Y}_n$ and $\widetilde{N}_n:=N_n/n$ by using the $ODE$ method of $SA$ under slightly lighten assumption than in \cite{BaiHu2}. The rate of convergence is investigated in Section~\ref{deux3}: we obtain a $CLT$, once again under slightly less stringent assumptions on the limit generating matrix $H$ than in \cite{BaiHu2}. Section~\ref{trois} is devoted to multi-arm clinical tests. In Section~\ref{trois1} we briefly recall the Wei $GFU$ model introduced~\cite{Wei,BaiHuShe} where the generating matrices $H_n$ are not random. In this case, the strong consistency and the asymptotic normality follow from the results of Section~\ref{deux} (like in \cite{BaiHu2}). In Section~\ref{trois2} we study the adaptive design proposed in~\cite{BaiHuShe} where the addition rule matrices depend on the responses of all the past patients. We use the results from Section~\ref{deux2} to prove the strong consistency. We prove in Section~\ref{trois3} a new $CLT$ for this model, when the generating matrix $H_n$ satisfies itself a $CLT$, which relies again on Stochastic Approximation techniques.

\medskip
\noindent  {\sc Notations} $\forall \, u=(u^i)_{i=1,\ldots,d}\in\R^d$, $\left\| u\right\|$ denotes the canonical Euclidean norm of the column vector $u$ on $\R^d$, $\Tr(u)=\sum_{k=1}^du^k$ denotes  its ``weight'',  $u^t$ denotes its transpose; $|\!|\!|A|\!|\!|$ denotes the operator norm of the matrix $A\in{\cal M}_{d,q}(\R)$ with $d$ rows and $q$ columns with respect to canonical Euclidean norms. When $d\!=\!q$, ${\rm Sp}(A)$ denotes the set of eigenvalues of $A$. $\mathbf{1}\!=\! (1\cdots1)^t$ denotes the unit column vector in $\R^d$, $I_d$ denotes the $d\times d$ identity matrix  and ${\rm diag}(u)=[\delta_{ij}u_i]_{1\le i,j\le d}$, where $\delta_{ij}$ is the Kronecker symbol. ${\cal S}=\left\{u\in\R^d_+:\sum_{i=1}^du^i=1\right\}$ denotes the $d$-dimensional simplex and ${\cal V}_0=\left\{u\in\R^d:\sum_{i=1}^du^i=0\right\}$.
 
\section{Convergence and first rate result}
\label{deux}

With the notations and definitions described in the introduction, we then formulate the main assumptions to establish the $a.s.$ convergence of the urn composition. 

\medskip
\noindent ${\bf (A1)}$ 
$\equiv\left\{\begin{array}{ll}
	(i) &  \mbox{{\em Addition rule matrix:}  For every $n\ge 1$, the  matrix $D_n$ $a.s.$ has non-negative entries.}\\
	\\
(ii) &  \mbox{{\em Generating matrix:}  For every $n\ge 1$,  the generating matrices $H_n=(H^{ij}_n)_{1\leq i,j \leq d}$} \\
	    &   \mbox{$a.s.$ satisfies} \\
	    &  \hskip2cm\forall\,j \in\{1,\ldots,d\}, \quad \displaystyle\sum_{i=1}^d H^{ij}_n=c>0. \\
	   
	(iii) & \mbox{{\em Starting value:} The starting urn composition vector $Y_0\!\in \R_+^d\setminus\{0\}$.}
\end{array}\right.$

\bigskip
The constant $c$ is known as the balance of the urn. In fact, we may assume without loss of generality, up to a renormalization of $Y_n$, that $c=1$: since $\widehat{Y}_n=\frac{Y_n}{c}$ and $\widehat{D}_{n+1}=\frac{D_{n+1}}{c}$, $n\geq 0$, formally satisfies the dynamics~(\ref{dynamic}), namely
$$
\widehat{Y}_n=\widehat{Y}_{n-1}+\widehat{D}_nX_n,\quad n\geq1, \quad \widehat{Y}_0\!\in \R_+^d\setminus\{0\}.
$$
From now on, throughout the paper, we will considered this normalized balance version. Nevertheless, we will still denote by $Y_n$ and $D_n$ the normalized quantities and assume that $c=1$.

\medskip
\noindent ${\bf (A2)}$ The addition rule $D_n$ is conditionally independent of the drawing procedure $X_n$ given $\F_{n-1}$ and satisfies
\begin{equation}\label{A2}
	\forall 1\leq j\leq d, \quad\sup_{n\geq1}\E\left[\left\|D^{\cdot j}_n\right\|^{2}\,|\,\F_{n-1}\right]<+\infty\quad a.s.
\end{equation}
where $D^{\cdot\,j}_n=(D^{ij}_n)_{i=1,\ldots,d}$.

The conditional independence is obtained in practice by assuming that the sequences of addition rules $(D_n)_{n\geq1}$ and the sequence $(U_n)_{n\geq1}$ used to randomize the drawings in (\ref{ConstructX}) are independent. \\

\noindent ${\bf (A3)}$ Assume that there exists an irreducible $d\times d$ matrix $H$ (with non-negative entries) such that 
\begin{equation}\label{A3}
	H_n\overset{a.s.}{\underset{n\rightarrow+\infty}{\longrightarrow}} H.
\end{equation}
$H$ is called the {\em limit generating matrix}. 

\medskip
The combination of assumptions {\bf (A1)}-{\bf (A3)} guarantees that $H$ satisfies the assumptions of the Perron-Frobenius Theorem (see \cite{BerPle}) so that $1$ is the eigenvalue of $H$ with the highest norm (maximal eigenvalue) has order $1$,  the components of its right eigenvector $v$ can be chosen  all positive and all other eigenvalues has a modulus  lower than $1$. In particular, we may normalize this vector $v^*$ such that $\Tr(v^*)=1$.

\bigskip
\noindent{\sc A variant including possible definite removal.} We may relax Assumption {\bf (A1)} by allowing the   removal of the drawn ball from its urn (see $e.g.$~\cite{Jan}). Other relaxation of these requirements may be considered: it could be possible to remove other balls than the drawn one. This leads to {\em tenable} urns  (studied notably in \cite{BagPal}, see also \cite{Pou}) where an arithmetical assumption to the row of any negative diagonal entry in $D_n$ is added, in order to avoid the urn extinction (see Assumption ${\bf (A'1)}$ below). Thus we may replace Assumption {\bf (A1)} (after renormalization) by\\

${\bf (A'1)}$
$\equiv\left\{\begin{array}{ll}
	(i) &  \mbox{{\em Addition rule matrix:} For every $i\!\in \{1,\ldots,d\}$, there exists $c_i\!\in (0,+\infty)$ such that,}\\
	    & \mbox{for every $n\ge 1$,}\\
	    & \forall\,i,\, j \in\{1,\ldots,d\}, \,\, \displaystyle  \frac{\delta_{ij}}{c_i} +D_n^{ij} \!\in \frac{\N}{c_i}\,\, a.s. \,\, \mbox{ and }\,\, \forall \, j \in\{1,\ldots,d\}, \,\, \sum_{i=1}^d D^{ij}_n \geq 0\quad a.s. \\
	(ii)&  \mbox{{\em Generating matrix:}  For every $n\ge 1$,   $H_n$ $a.s.$ satisfies} \\
	    
	    & \hskip3cm\forall\,j \in\{1,\ldots,d\}, \quad \displaystyle\sum_{i=1}^d H^{ij}_n=1. \\	    
	(iii) & \mbox{{\em Starting value:} The starting urn composition vector $\displaystyle Y_0\!\in\Big( \prod_{i=1}^d\frac{\N}{c_i}\Big)\setminus\{0\}$.} \\
\end{array}\right.$
 
\bigskip
In this case $H$ may have negative (diagonal) entries and the Perron-Frobenius Theorem cannot be used, so we change Assumption {\bf (A3)} into 
\begin{equation*}
{\bf (A'3)} \; \mbox{ $1$ is the eigenvalue of $H$ with maximal modulus, has order $1$ and  $\{v:Hv=v\}\subset \R_+^d$.\hskip3cm }
\end{equation*}
Throughout the paper, we may substitute ${\bf (A'1)}$-${\bf (A'3)}$ for {\bf (A1)}-{\bf (A3)} as recalled in each result.
 
The following preliminary lemma ensures that if ${\bf (A'1)}$ holds then the urn extinction never occurs and its weight $\Tr (Y_n)$ is non-decreasing.
\begin{lem}[Preliminary] If ${\bf (A'1)}$ holds, then $\Tr(Y_n)$ is non-decreasing and postive.
\end{lem}

\noindent{\bf Proof.} We proceed by induction on $n\ge 0$. Assume $\displaystyle Y_{n-1} \!\in\Big( \prod_{i=1}^d \frac{\N}{c_i}\Big)\setminus\{0\}$. For every $i\!\in \{1,\ldots,d\}$, 
 \[
 Y_n^i =Y^i_{n-1} +\sum_{j=1}^dD_n^{ij}\mathds{1}_{\{X_n=e^j\}} \quad \mbox{and} \quad\{X_n=e^j\}\subset \{Y^j_{n-1}>0\}= \{Y^j_{n-1} \ge 1/c_j\}.
 \] 
 Consequently $Y^i_n \ge Y^i_{n-1}$ and $\displaystyle Y^i_n \!\in \frac{\N}{c_i}\setminus\{0\}$ on the event $\bigcup_{j\neq i} \{X_n=e^j\}$. On $\{X_n=e^i\}$, $\{Y^i_{n-1}\ge \frac{1}{c_i}$ so that $Y^i_n=Y^i_{n-1}+D_n^{ii} \ge \frac{1}{c_i}-\frac{1}{c_i}\ge 0$. 
Finally
 \[
 \Tr(Y_n)=\Tr(Y_{n-1}) + \sum_{j=1}^d \Big(\sum_{i=1}^d D_n^{ij}\Big)\mathds{1}_{\{X_n=e^j\}}\ge \Tr (Y_{n-1})>0.\qquad \hfill \cqfd
 \] 

\subsection{The dynamics as a stochastic approximation procedure}
\label{deux1}

Our aim in this section is to reformulate the dynamics (\ref{dynamic})-(\ref{ConstructX}) into a recursive stochastic algorithm. Then we aim at applying the most powerful tools of $SA$, namely the ``$ODE$"   and the ``$SDE$"  methods to elucidate the asymptotic properties ($a.s.$ convergence and weak rate) of both the urn composition and the treatment allocation. We start from~(\ref{dynamic}) with $Y_0\!\in\R_+^d\setminus\{0\}$. For $n\geq1$,
\begin{equation}\label{dynq}
Y_{n+1}=Y_n+D_{n+1}X_{n+1}=Y_n+\E\left[D_{n+1}X_{n+1}\,|\,\F_{n}\right]+\Delta M_{n+1},
\end{equation}
where 
$$
\Delta M_{n+1}:=D_{n+1}X_{n+1}-\E\left[D_{n+1}X_{n+1}\,|\,\F_{n}\right]
$$ 
is an $\F_{n}$-martingale increment. By the definition of the generating matrix $H_n$, we have
\begin{eqnarray*}
\E\left[D_{n+1}X_{n+1}\,|\,\F_{n}\right]&=&\sum_{i=1}^d\E\left[D_{n+1}\mathds{1}_{\left\{X_{n+1}=e^i\right\}}e^i\,|\,\F_{n}\right]=\sum_{i=1}^d\E\left[D_{n+1}\,|\,\F_{n}\right]\P\left(X_{n+1}=e^i\,|\,\F_{n}\right)e^i\\
 &=&H_{n+1}\sum_{i=1}^d\frac{Y_n^i}{\Tr(Y_n)}e^i=H_{n+1}\frac{Y_n}{\Tr(Y_n)}
\end{eqnarray*}
so that $\displaystyle \hskip 2cm Y_{n+1}=Y_n+H_{n+1}\frac{Y_n}{\Tr(Y_n)}+\Delta M_{n+1}$. 

Now we can derive a stochastic approximation for the normalized urn composition $Y_n$. First we have for every $n\ge 1$, 
$$
\frac{Y_{n+1}}{n+1}=\frac{Y_n}{n}+\frac{1}{n+1}\left(H_{n+1}\frac{Y_n}{\Tr(Y_n)}-\frac{Y_n}{n}\right)+\frac{\Delta M_{n+1}}{n+1}.
$$
Consequently,  $\displaystyle \widetilde{Y}_n=\frac{Y_n}{n}$, $n\ge 1$, satisfies a canonical recursive stochastic approximation procedure
\begin{eqnarray}
\widetilde{Y}_{n+1}&=&\widetilde{Y}_n+\frac{1}{n+1}\left(H_{n+1}-I_d\right)\widetilde{Y}_n+\frac{1}{n+1}\left(\Delta M_{n+1}+\left(\frac{n}{\Tr(Y_n)}-1\right)H_{n+1}\widetilde{Y}_n\right)\nonumber\\
&=&\widetilde{Y}_n-\frac{1}{n+1}\left(I_d-H\right)\widetilde{Y}_n+\frac{1}{n+1}\left(\Delta M_{n+1}+r_{n+1}\right) \label{AS}
\end{eqnarray}	                   
with step $\gamma_n =\frac 1n$ and a remainder term given by
\begin{equation}\label{reste}
	r_{n+1}:=\left(\frac{n}{\Tr(Y_n)}-1\right)H_{n+1}\widetilde{Y}_n+(H_{n+1}-H)\widetilde{Y}_n. 
\end{equation}

Furthermore, in order to establish the $a.s.$ boundedness of $(\widetilde{Y}_n)_{n\geq1}$ we will rely on the following recursive equation satisfied by $\Tr(Y_n)$:
$$\Tr(Y_{n+1})=\Tr(Y_n)+\frac{\Tr(H_{n+1}Y_n)}{\Tr(Y_n)}+\Tr(\Delta M_{n+1}).$$
By the properties of the generating matrix $H_{n+1}$, we obtain
$$\Tr(H_{n+1}Y_n)=\sum_{i=1}^d(H_{n+1}Y_n)_i=\sum_{i=1}^d\sum_{j=1}^dH^{ij}_{n+1}Y_n^j=\sum_{j=1}^d\left(\sum_{i=1}^dH^{ij}_{n+1}\right)Y_n^j=\Tr(Y_n).$$
Consequently 
\begin{equation}\label{ASTrace}
\Tr(Y_{n+1})=\Tr(Y_n)+1+\Tr(\Delta M_{n+1}).
\end{equation}

\subsection{Convergence results}
\label{deux2}

\begin{theo}\label{Thm1}
Let $(Y_n)_{n\geq0}$ be the urn composition sequence defined by~(\ref{dynamic})-(\ref{ConstructX}). Under the assumptions {\bf (A1)}, {\bf (A2)} and {\bf (A3)} (or ${\bf (A'1)}$, {\bf (A2)} and ${\bf (A'3)}$),
\begin{enumerate}
	\item[$(a)$] $\frac{\Tr(Y_n)}{n}\overset{a.s.}{\underset{n\rightarrow+\infty}{\longrightarrow}} 1 $ and $\displaystyle\frac{Y_n}{\Tr(Y_n)}\overset{a.s.}{\underset{n\rightarrow+\infty}{\longrightarrow}}v^*$.
	\item[$(b)$] $\widetilde{N}_n:=\displaystyle \frac{N_n}{n}=\displaystyle\frac{1}{n}\sum_{k=1}^nX_k\overset{a.s.}{\underset{n\rightarrow+\infty}{\longrightarrow}}v^*$.
\end{enumerate}
\end{theo}

\medskip \noindent {\bf Remarks.} $\bullet$ We simply need that $H_n\overset{a.s.}{\underset{n\rightarrow+\infty}{\longrightarrow}}H$ while the assumption in~\cite{BaiHu2} is
$$
\sum_{n\geq1}\frac{\left\|H_n-H\right\|_{\infty}}{n}<+\infty 
$$
where $\left\|\cdot\right\|_{\infty}$ is the norm on $L^{\infty}_{\R^{d\times d}}(\P)$. \\
\noindent $\bullet$ Assumption {\bf (A3)} is not necessary to prove that $\frac{\Tr(Y_n)}{n}\overset{a.s.}{\underset{n\rightarrow+\infty}{\longrightarrow}}1$.

\bigskip
\noindent {\bf Proof.}  We will first prove that $(a)\Rightarrow(b)$, then we will prove $(a)$. 

\smallskip
\noindent $(a)\Rightarrow(b).$ We have
$$
\E\left[X_n\,|\,\F_{n-1}\right]=\sum_{i=1}^d\frac{Y_{n-1}^i}{\Tr(Y_{n-1})}e^i=\frac{Y_{n-1}}{\Tr(Y_{n-1})}
$$ 
and, by construction $\left\|X_n\right\|^2=1$ so that $\E\left[\left\|X_n\right\|^2\,|\,\F_{n-1}\right]=1$. Hence the martingale
$$
\widetilde{M}_n=\sum_{k=1}^n\frac{X_k-\E\left[X_k\,|\,\F_{k-1}\right]}{k}\overset{a.s.\&\, L^2}{\underset{n\rightarrow+\infty}{\longrightarrow}}\widetilde{M}_{\infty}\in L^2,
$$
and by the Kronecker Lemma we obtain
$$
\frac{1}{n}\sum_{k=1}^nX_k-\frac{1}{n}\sum_{k=1}^n\frac{Y_{k-1}}{\Tr(Y_{k-1})}\overset{a.s.}{\underset{n\rightarrow+\infty}{\longrightarrow}}0.
$$
This yields the announced implication owing  to the Cesaro Lemma. 

\medskip
\noindent $(a)$ {\sc First Step}: We have 
$$
D_{n+1}X_{n+1}=\sum_{j=1}^dD^{\cdot\,j}_{n+1}\mathds{1}_{\{X_{n+1}=e^j\}}.
$$
Therefore
\begin{eqnarray*}
\left\|D_{n+1}X_{n+1}\right\|^2&=&\sum_{j=1}^d\left\|D^{\cdot\, j}_{n+1}\right\|^2\mathds{1}_{\{X_{n+1}=e^j\}},\\
\mbox{so that } \hskip 1 cm 	\E\left[\left\|D_{n+1}X_{n+1}\right\|^2\,|\,\F_n\right]&=&\sum_{j=1}^d\E\left[\left\|D^{\cdot j}_{n+1}\right\|^2\,|\,\F_n\right]\P\left(X_{n+1}=e^j\,|\,\F_n\right) \\
	&\leq&\,\sup_{n\geq0}\sup_{1\leq j\leq d}\E\left[\left\|D^{\cdot j}_{n+1}\right\|^2\,|\,\F_n\right]<+\infty\quad a.s.\hskip 2 cm 
\end{eqnarray*}
Consequently $\sup_{n\geq1}\E\left[\left\|\Delta M_{n+1}\right\|^2\,|\,\F_n\right]<+\infty$ $a.s.$. 
Therefore thanks to the strong law of large numbers for conditionally $L^2$-bounded martingale increments, we have $\frac{M_n}{n}\overset{a.s.}{\underset{n\rightarrow+\infty}{\longrightarrow}}0$. Consequently it follows from (\ref{ASTrace}) that
\begin{equation}\label{ConvTr}
\frac{\Tr(Y_n)}{n}=1+\frac{\Tr(Y_0)-1}{n}+\frac{\Tr(M_n)}{n}\overset{a.s.}{\underset{n\rightarrow+\infty}{\longrightarrow}} 1.
\end{equation}

\medskip
\noindent {\sc Second Step:} Since the components of $\widetilde{Y}_n=\frac{Y_n}{n}$ are non-negative and $\Tr(\widetilde{Y}_n)=\frac{\Tr(Y_n)}{n}\overset{a.s.}{\underset{n\rightarrow+\infty}{\longrightarrow}} 1$, it is clear that $(\widetilde{Y}_n)_{n\geq1}$ is $a.s.$ bounded and that $a.s.$ the set ${\cal Y}_{\infty}$ of all its limiting value is contained in
$${\cal S}=\Tr^{-1}\{1\}=\left\{u\in\R^d_+\,|\,\Tr(u)=1\right\}.$$
So we may try applying the $ODE$ method (see Appendix Theorem~\ref{ThmODE}). Since $\widetilde{Y}_n$ and $H_{n+1}\widetilde{Y}_n$ are $a.s.$ bounded, (\ref{ConvTr}) and {\bf (A3)} imply that $r_n\overset{a.s.}{\underset{n\rightarrow+\infty}{\longrightarrow}}0$.

The $ODE$ associated to the recursive procedure reads
$$ODE_{I_d-H}\,\equiv\,\dot{y}=-(I_d-H)y.$$
Owing to Assumption {\bf (A3)}, $I_d-H$ admits $v^*$ as unique zero in ${\cal S}$. The restriction of $ODE_{I_d-H}$ to the affine hyperplane ${\cal V}$ is the linear system $\dot{z}=-(I_d-H)z$, where $z=y-v^*$ takes values in ${\cal V}_0=\left\{u\in\R^d\,|\,\Tr(u)=0\right\}$. Since $\mbox{Sp}\left((I_d-H)_{\,|\,{\cal V}_0}\right)\subset\left\{\lambda\in\C, \, \Re e(\lambda)>0\right\}$, owing to Assumption~{\bf (A3)}. As a consequence $v^*$ is a uniformly stable equilibrium for the restriction of $ODE_{I_d-H}$ to ${\cal S}$, the whole hyperplane, as an attracting area. The fundamental result derived from the $ODE$ method (see Theorem~\ref{ThmODE} in Appendix and the notations therein, in particular the remainder $r_n$) yields the expected result
$$
\hskip7cm\widetilde{Y}_n\overset{a.s.}{\underset{n\rightarrow+\infty}{\longrightarrow}}v^*.\hskip6.5cm\cqfd
$$

\paragraph{Remark:} If we assume that the addition rule matrices $(D_n)_{n\ge 1}$ satisfy besides {\bf (A1)}, then we can directly write a stochastic approximation for $\frac{Y_n}{\Tr(Y_n)}$ with step $\frac{1}{\Tr(Y_n)}$ in which the remainder simply reads $r_{n+1}=(H_{n+1}-H)\frac{Y_n}{\Tr(Y_n)}$ and prove the $a.s.$ convergence under the same assumptions.

\medskip
\noindent {\sc Comments.} We could apply directly the $ODE$ method because we first proved that $(\widetilde{Y}_n)_{n\geq1}$ is $a.s.$ bounded without using the standard Lyapunov machinery developed in $SA$ Theory. That is why the assumption on the remainder sequence $(r_n)_{n\geq1}$ simply reads 
$$r_n\overset{a.s.}{\underset{n\rightarrow+\infty}{\longrightarrow}}0.$$
Another approach is the martingale one. It relies on the existence of a Lyapunov function $V:\R^d\rightarrow\R_+$ associated to the algorithm satisfying
\begin{equation}\label{Lyap}
\exists\, a>0,\quad\forall y\in\R^d,\quad y\neq v^*,\quad \left\langle \nabla V\left|\right.I_d-H\right\rangle(y)>0 \quad\mbox{and}\quad\left\langle \nabla V\left|\right.I_d-H\right\rangle>a\left|\nabla V\right|^2.
\end{equation}
In this framework the existence of a Lyapunov function can be established. Hence, the natural condition on the remainder sequence $(r_n)_{n\geq1}$ reads (see \cite{Duf})
$$
\sum_{n\geq1}\frac{\left\|r_n\right\|^2}{n}<+\infty\quad a.s.
$$
In that perspective, the assumption on the generating matrices would read $\displaystyle \sum_{n\geq1}\frac{|\!|\!|H_n-H|\!|\!|^2}{n}<+\infty$ $a.s.$
which is still slightly less stringent than assumption on the generating matrices made in~\cite{BaiHu2}.

\subsection{Rate of convergence}
\label{deux3}

In the previous section we proved the $a.s.$ convergence of both quantities of interest, namely $\widetilde{Y}_n$ and $\widetilde{N}_n$, toward $v^*$. In this section we establish a ``joint $CLT$'' for the (column) couple 
$$
\theta_n:=(\widetilde{Y}_n, \widetilde{N}_n)^t
$$ 
with an explicit  asymptotic joint normal distribution (including covariances). To this end we will show that $\theta_n$ satisfies a ${\cal S}^2$-valued $SA$ recursive procedure which ($a.s.$ converges toward $\theta^*=\left(v^*, v^*\right)^t\in{\cal S}^2$ and) fulfills the assumptions of the $CLT$ Theorem~\ref{ThmCLT} for $SA$ algorithms (see Appendix A ans Appendix C for the spectrum of $Dh(\theta^*)$), with a special attention paid to Condition~(\ref{HypReste}) about the  remainder  term. 
As concerns $\widetilde Y_n$, we derive from (\ref{AS}) that
$$
\forall\, n\ge 1, \qquad \widetilde{Y}_{n+1}=\widetilde{Y}_n-\frac{1}{n+1}\left(I_d-(2-\Tr(\widetilde{Y}_n))H\right)\widetilde{Y}_n+\frac{1}{n+1}\left(\Delta M_{n+1}+\bar{r}_{n+1}\right),
$$	
where $\hskip3.5cm \bar{r}_{n+1}:=\displaystyle\left(\frac{H_{n+1}-H}{\Tr(\widetilde{Y}_n)}+\frac{(\Tr(\widetilde{Y}_n)-1)^2}{\Tr(\widetilde{Y}_n)}H\right)\widetilde{Y}_n$. 

\noindent For $\widetilde{N}_n$ we have, still for every $n\ge 1$, 
\begin{equation}\label{eq:Ntilde}
\widetilde{N}_{n+1}=\widetilde{N}_n-\frac{1}{n+1}\left(\widetilde{N}_n-(2-\Tr(\widetilde{Y}_n))\widetilde{Y}_n\right)+\frac{1}{n+1}\left(\Delta \widetilde{M}_{n+1}+\widetilde{r}_{n+1}\right)
\end{equation}
with $\displaystyle\Delta \widetilde{M}_{n+1}:=X_{n+1}-\E\left[X_{n+1}\left|\right.\F_n\right]=X_{n+1}-\frac{Y_n}{\Tr(Y_n)}$ and $\widetilde{r}_{n+1}:=\displaystyle\frac{(\Tr(\widetilde{Y}_n)-1)^2}{\Tr(\widetilde{Y}_n)}\widetilde{Y}_n$.

\noindent Thus, we obtain a new recursive $SA$ procedure, still with step $\gamma_n = \frac 1n$, namely
$$
\theta_{n+1}=\theta_n-\frac{1}{n+1}h(\theta_n)+\frac{1}{n+1}\left(\Delta {\bf M}_{n+1}+R_{n+1}\right), \quad n\ge1,
$$
with $\Delta {\bf M}_{n+1}:=\begin{pmatrix} \Delta M_{n+1} \cr \Delta \widetilde{M}_{n+1}\end{pmatrix}$, $R_{n+1}:=\begin{pmatrix} \bar{r}_{n+1} \cr  \widetilde{r}_{n+1}\end{pmatrix}$ and
$$
\forall\, \theta=\begin{pmatrix} y \cr \nu\end{pmatrix},\,y\in\R^d,\,\nu\in\R^d, \quad h(\theta):=\begin{pmatrix} (I_d-(2-\Tr(y))H)y \cr \nu-(2-\Tr(y))y \end{pmatrix} \quad \mbox{with} \quad h(\theta^*)=0.
$$

\noindent The function $h$ is differentiable on $\R^{d}\times \R^d$ and its differential at point $\theta^*$  is given by
$$
Dh(\theta^*)=\begin{pmatrix} I_d-H+v^*\mathbf{1}^t &   0_{{\cal M}_{d}(\R)}\cr   v^*\mathbf{1}^t-I_d &  I_d\end{pmatrix} \; \mbox{so that }\; Dh(\theta^*)_{|{\cal V}_0^2}=\begin{pmatrix} (I_d-H)_{|\mbox{\bf  1}^\perp} &  0_{_{|\mbox{\bf  1}^\perp} }    \cr  - I_{d|\mbox{\bf  1}^\perp}&  I_{d|\mbox{\bf  1}^\perp} \end{pmatrix}.
$$

\noindent To establish a $CLT$ for the sequence $(\theta_n)_{n\geq 1}$ we need to make  the following additional assumptions:

\medskip
\noindent {\bf (A4)} The addition rules $D_n$ $a.s.$ satisfy
\begin{equation*}\label{A4}
	\forall 1\leq j\leq d, \quad\left\{
	\begin{array}{ll}
	 	\sup_{n\geq1}\E\left[\|D^{\cdot j}_n\|^{2+\delta}\,|\,\F_{n-1}\right]\leq C<+\infty & \mbox{for a  $\delta>0$,} \\
		\E\left[D^{\cdot j}_n(D^{\cdot j}_n)^t\,|\,\F_{n-1}\right]\underset{n\rightarrow+\infty}{\longrightarrow} C^j, & \\
	\end{array}\right.
\end{equation*}
where $C^j=(C_{il}^j)_{1\leq i,l\leq d}$, $j=1,\ldots,d$, are $d\times d$ symmetric  positive definite matrices. 

Note that {\bf (A4)} $\Rightarrow${\bf (A2)} since $\E\left[\|D^{\cdot j}_n\|^2\,|\,\F_{n-1}\right]\leq\left(\E\left[\|D^{\cdot j}_n\|^{2+\delta}\,|\,\F_{n-1}\right]\right)^{\frac{2}{2+\delta}}$.

\medskip
\noindent {\bf (A5)$_v$} The matrix $H$ satisfies
\begin{equation}\label{A5}
	 nv_n\,\E\left[|\!|\!|H_n-H|\!|\!|^2\right]\underset{n\rightarrow+\infty}{\longrightarrow}0,
\end{equation}
where $(v_n)_{n\geq1}$ is a positive  sequence (specified in each item of the theorems further on).

\begin{theo}\label{Thm2}
Assume {\bf (A1)}, {\bf (A3)} (or ${\bf (A'1)}$, ${\bf (A'3)}$), {\bf (A4)} and {\bf (A5)}. 

\noindent$(a)$ Let $\lambda_{\max}$  the eigenvalue of H with the highest real part appart from $1$. If
\begin{equation}\label{stable}
\lambda_{\max} = \max\Re e\left({\rm Sp}(H)\setminus\{1\}\right) <1/2
\end{equation}
and~\eqref{A5} holds with  $v_n=1$, $n\geq1$, then, $\theta_n \to \theta^*$ a.s. and 
$$
\sqrt{n}\left(\theta_n-\theta^*\right)\overset{\L}{\longrightarrow}{\cal N}\left(0,\Sigma\right)\, \mbox{as $n\to+\infty\,$ with }\Sigma=\int_0^{+\infty}e^{-u\left(Dh(\theta^*)-\frac{I_{2d}}{2}\right)^t}\Gamma e^{-u\left(Dh(\theta^*)-\frac{I_{2d}}{2}\right)}du
$$
 \begin{equation}\label{Gamma}
\mbox{and }\;\Gamma=\!\begin{pmatrix}\displaystyle\sum_{k=1}^d v^{*k}C^k-v^*(v^*)^t & H\left({\rm diag}(v^*)-v^*(v^*)^t\right) \cr & & \cr \left({\rm diag}(v^*)-v^*(v^*)^t\right)^tH^t & {\rm diag}(v^*)-v^*(v^*)^t\!\!\end{pmatrix}=a.s.\mbox{-}\lim_{n\rightarrow+\infty}\E\left[\Delta {\bf M}_n\Delta {\bf M}_n^t\,|\,\F_{n-1}\right].
\end{equation}
$(b)$ If $\lambda_{{\rm max}}=1/2$, $H$ is $\R$-diagonalizable and~\eqref{A5} holds with  $v_n=\log n$, $n\geq 2$, then $\theta_n \to \theta^*$ $a.s.$  and 
$$
\sqrt{\frac{n}{\log n}}\left(\theta_n-\theta^*\right)\overset{\L}{\underset{n\rightarrow+\infty}{\longrightarrow}}{\cal N}\left(0,\Sigma\right)\quad\mbox{ with }\quad \Sigma=\lim_{n\to+\infty}\frac{1}{\log n}\int_0^{\log n}e^{-u\left(Dh(\theta^*)-\frac{I_{2d}}{2}\right)^t}\Gamma e^{-u\left(Dh(\theta^*)-\frac{I_{2d}}{2}\right)}du.
$$
$(c)$ If $ \lambda_{{\rm max}}\!\in(1/2,1)$, $H$ is $\R$-diagonalizable and~\eqref{A5} holds with $v_n=n^{1-2\lambda_{\max}+\eta}$, $n\geq1$, for some $\eta>0$, then $\theta_n \to \theta^*$ a.s. and $n^{1-\lambda_{{\rm max}}}\left(\theta_n-\theta^*\right)$ $a.s.$ converges as $n\to+\infty$ towards a finite random variable.
\end{theo}

\noindent {\bf Proof.} $(a)$ We will check the three assumptions of the $CLT$  for $SA$ algorithms recalled in the Appendix (Theorem~\ref{ThmCLT}). Firstly, the condition~(\ref{HypoPas}) on the spectrum of $Dh(\theta^*)_{|{\cal S}} $ requested for algorithms with step $\frac 1n$ in Theorem~\ref{ThmCLT} reads $\Re e\left({\rm Sp}(Dh(\theta^*)_{|{\cal S}} )\right)>\frac 12$. This follows from our Assumption~(\ref{stable}) since by decomposing $\R^d= \R v^*\oplus {\rm Ker}(\Tr)$, one checks that 
\[
{\rm Sp}(Dh(\theta^*)_{|{\cal S}} )= \{1\}\cup\big\{1-\lambda,\, \lambda\!\in {\rm Sp}(H)\setminus\{1\}\big\}.
\] 

\noindent Secondly Assumption~{\bf (A4)} ensures that Condition~(\ref{HypDM}) is satisfied since
$$
\sup_{n\geq1}\E\left[\left\|\Delta {\bf M}_n\right\|^{2+\delta}\,|\,\F_{n-1}\right]<+\infty \quad a.s.\quad \mbox{and} \quad
\E\left[\Delta {\bf M}_n\Delta {\bf M}_n^t\left.\right|\F_{n-1}\right]\overset{a.s.}{\underset{n\rightarrow+\infty}{\longrightarrow}}\Gamma\quad \mbox{as}\quad n\to +\infty,
$$
where $\Gamma$ is the symmetric nonnegative matrix given by (\ref{Gamma}) as established below. To this end we have to determine three blocks since $\Gamma$ reads
$$
\Gamma=\begin{pmatrix}\Gamma_1 & \Gamma_{12} \cr \Gamma_{12}^t & \Gamma_{2}\end{pmatrix}\quad\mbox{where}\quad \Gamma_1,\Gamma_2,\Gamma_{12}\in{\cal M}_{d}(\R).
$$
\noindent {\em Computation of $\Gamma_1$.}

\vskip -0.75cm
\begin{eqnarray*}
\E\left[\Delta M_{n+1}\Delta M_{n+1}^t\,|\,\F_n\right] \!\!&\!\!=\!\!&\!\! \sum_{q=1}^d\P(X_{n+1}=e^q\,|\,\F_n) \! \left(\E\left[D_{n+1}^{\cdot q}(D_{n+1}^{\cdot q})^t\,|\,\F_n\right]\!\right. \\
                                                           & &-\left.\!\E\left[D_{n+1}X_{n+1}\,|\,\F_n\right]\E\left[D_{n+1}X_{n+1}\,|\,\F_n\right]^t\right) \\
	 \!\!&\!\!=\!\!&\!\!\sum_{q=1}^d\frac{Y_n^q}{\Tr(Y_n)}\E\left(D_{n+1}^{\cdot q}(D_{n+1}^{\cdot q})^t\,|\,\F_n\right)-\left(H_{n+1}\frac{Y_n}{\Tr(Y_n)}\right)\left(H_{n+1}\frac{Y_n}{\Tr(Y_n)}\right)^t \\
	 &\overset{a.s.}{\underset{n\rightarrow+\infty}{\longrightarrow}}&\Gamma_1=\sum_{q=1}^d v^{*q}C^q-v^*(v^*)^t.
\end{eqnarray*}

\noindent {\em Computation of $\Gamma_2$.}
\begin{eqnarray*}
	\E\left[\Delta \widetilde{M}_{n+1}\Delta \widetilde{M}_{n+1}^t\,|\,\F_n\right]&=&\E\left[X_{n+1}X_{n+1}^t\,|\,\F_n\right]-\frac{Y_n}{\Tr(Y_n)}\left(\frac{Y_n}{\Tr(Y_n)}\right)^t \\
	&=&\mbox{diag}\left(\frac{Y_n}{\Tr(Y_n)}\right)-\frac{Y_n}{\Tr(Y_n)}\left(\frac{Y_n^q}{\Tr(Y_n)}\right)^t\overset{a.s.}{\underset{n\rightarrow+\infty}{\longrightarrow}}\Gamma_2=\mbox{diag}(v^*)-v^*(v^*)^t.
\end{eqnarray*}

\noindent {\em Computation of $\Gamma_{12}$.}
\begin{eqnarray*}
		\E\left[\Delta M_{n+1}\Delta \widetilde{M}_{n+1}^t\,|\,\F_n\right]&=&\E\left[D_{n+1}X_{n+1}X_{n+1}^t\,|\,\F_n\right]-\E\left[D_{n+1}X_{n+1}\,|\,\F_n\right]\E\left[X_{n+1}\,|\,\F_n\right]^t \\	&=&\E\left[D_{n+1}\,|\,\F_n\right]\E\left[X_{n+1}X_{n+1}^t\,|\,\F_n\right]-\E\left[D_{n+1}\,|\,\F_n\right]\E\left[X_{n+1}\,|\,\F_n\right]\E\left[X_{n+1}\,|\,\F_n\right]^t \\	&=&H_{n+1}\mbox{diag}\left(\frac{Y_n}{\Tr(Y_n)}\right)-H_{n+1}\frac{Y_n}{\Tr(Y_n)}\left(\frac{Y_n}{\Tr(Y_n)}\right)^t\\
		&\overset{a.s.}{\underset{n\rightarrow+\infty}{\longrightarrow}}&\Gamma_{12}=H\left(\mbox{diag}(v^*)-v^*(v^*)^t\right).
\end{eqnarray*}

\noindent Finally, it remains to check that the remainder sequence $(R_n)_{n\geq1}$ satisfies~(\ref{HypReste}) for an $\epsilon>0$:
\begin{equation}\label{HReste}
\E\left[(n+1)\left\|R_{n+1}\right\|^2\mathds{1}_{\{\left\|\theta_n-\theta^*\right\|\leq \epsilon\}}\right]\underset{n\rightarrow+\infty}{\longrightarrow}0.
\end{equation}
We note that  $\left\|R_{n+1}\right\|^2=\left\|\bar{r}_{n+1}\right\|^2+\left\|\widetilde{r}_{n+1}\right\|^2$. It follows from the definition of $\bar{r}_{n+1}$ and the elementary facts  $\|\widetilde{Y}_n-v^*\|\leq \|\theta_n-\theta^* \|$ and $\Tr(\widetilde{Y}_n)\geq \|\widetilde{Y}_n \|$ that
\begin{eqnarray*}
	\left\|\bar{r}_{n+1}\right\|^2\mathds{1}_{\left\{\left\|\theta_n-\theta^*\right\|\leq \frac{\left\|v^*\right\|}{2}\right\}}&\leq&2\left(\frac{(\Tr(\widetilde{Y}_n)-1)^4}{\frac{\left\|v^*\right\|}{2}}+\frac{|\!|\!| H_{n+1}-H|\!|\!|^2}{\frac{\left\|v^*\right\|}{2}}\right)\frac 32\|v^*\|\mathds{1}_{\left\{\left\|\theta_n-\theta^*\right\|\leq \frac{\left\|v^*\right\|}{2}\right\}} \\
	&\leq&6\left((\Tr(\widetilde{Y}_n)-1)^4+|\!|\!|H_{n+1}-H|\!|\!|^2\right)\mathds{1}_{\left\{\left\|\theta_n-\theta^*\right\|\leq\frac{\left\|v^*\right\|}{2}\right\}}.
\end{eqnarray*}
But $\Tr(\widetilde{Y}_n)-1=\frac{\Tr(\Delta M_n)}{n}$ where $\sup_{n\geq0}\E\left[\left|\Tr(\Delta M_{n+1})\right|^{2+\delta}\left|\right.\F_n\right]\leq C'$, $\delta>0$, owing to (A4). Now using that $\left|\Tr(y)\right|\leq C_d\|y\|$,
$$
\E\left[n\left|\Tr(\widetilde{Y}_n)-1\right|^4\mathds{1}_{\left\{\left\|\theta_n-\theta^*\right\|\leq\frac{\left\|v^*\right\|}{2}\right\}}\right]\leq C^*_{\delta} n\E\left[\left|\Tr(\widetilde{Y}_n)-1\right|^{2+\delta}\right]=\frac{C_d}{n^{1+\delta}}\E\left[\left|\Tr(\Delta M_n)\right|^{2+\delta}\right]\leq\frac{C'_d}{n^{1+\delta}},$$
where $C^*_{\delta}>0$ is a real constant. Consequently 
$$
n\E\left[\left|\Tr(\widetilde{Y}_n)-1\right|^4\mathds{1}_{\left\{\left\|\theta_n-\theta^*\right\|\leq\frac{\left\|v^*\right\|}{2}\right\}}\right]=O\left(\frac{1}{n^{\delta}}\right).$$
Thus, by (A5) we obtain
$$
n\E\left[\left\|\bar{r}_{n+1}\right\|^2\mathds{1}_{\left\{\left\|\theta_n-\theta^*\right\|\leq\frac{\left\|v^*\right\|}{2}\right\}}\right]=O\left(\frac{1}{n^{\delta}}\right).
$$
The same argument yields  $ \displaystyle 
n\E\left[\left\|\widetilde{r}_{n+1}\right\|^2\mathds{1}_{\left\{\left\|\theta_n-\theta^*\right\|\leq\frac{\left\|v^*\right\|}{2}\right\}}\right]=O\left(\frac{1}{n^{\delta}}\right),
$
therefore the remainder condition~(\ref{HReste}) is satisfied. 

\smallskip
\noindent \noindent $(b)$-$(c)$ follow from Theorem~\ref{ThmCLT} $(b)$-$(c)$ in  the Appendix  since one easily checks that $D(h(\theta^*))_{|{\cal V}^2}$ is diagonalizable as soon as $H$ is with ${\rm Sp}( D(h(\theta^*))_{|{\cal S}^2})= \{1-\lambda,\, \lambda\!\in {\rm Sp}(H)\setminus \{1\}\}$ (see Appendix~\ref{app:C}). Moreover the above computations show that the remainder condition~(\ref{HReste}) is satisfied.~$\cqfd$

%

\section{Application to urn models for multi-arm clinical trials}
\label{trois}

In this section, we consider urn models for multi-arm clinical trials introduced by Wei and generalized by Bai, Hu and Shen. In this context, the initial framework where the addition rule matrices have nonnegative entries is the only one to make sense.

\subsection{The Wei {\em GFU} Model}
\label{trois1}
We consider here the model presented in~\cite{Wei} and in~\cite{BaiHuShe}, where balls are added depending on the success probabilities of each treatment. Define an {\em efficiency indicator} as follows: let $(T_n^i)_{n\geq1}$, $1\leq i\leq d$, be $d$ independent sequences of $[0,1]$-valued i.i.d. random variables, independent of the i.i.d.{\em sampling} sequence $(U_n)_{n\ge 1} $  so that 
\begin{equation}\label{ProbaSucces}
	\E\left[T^i_n\right]=p^i, \quad 0<p^i<1,\quad 1\leq i\leq d.
\end{equation}

\paragraph{Remark.} If $(T_n^i)_{n\geq1}$, $1\leq i\leq d$, is simply a {\em success indicator}, namely $d$ independent sequences of i.i.d. $\{0,1\}$-valued Bernoulli trials with respective parameter $p^i$, then the convention is to set $T_n^i\!=\!1$ to indicate that the response of the $i^{th}$ treatment in the $n^{th}$ trial is a success and $T_n^i=0$ otherwise. \\

In this framework one considers the filtration $\F_n=\sigma\left(Y_0,U_k,T_k,1\leq k\leq n\right)$, $n\ge 0$.
Consider the following addition rules: a success on the treatment $i$ adds a ball of type $i$ to the urn and a failure on the treatment $i$ adds $\frac{1}{d-1}$ balls for each of the other $d-1$ types. Thus the addition rule proposed in \cite{Wei} is as follows
$$
D_{n+1}=\begin{pmatrix}T^1_{n+1} & \frac{1-T^2_{n+1}}{d-1} & \cdots & \frac{1-T^d_{n+1}}{d-1} \cr \cr \frac{1-T^1_{n+1}}{d-1} & T^2_{n+1} & \cdots & \frac{1-T^d_{n+1}}{d-1} \cr \vdots & \vdots & \ddots & \vdots \cr \frac{1-T^1_{n+1}}{d-1} & \frac{1-T^2_{n+1}}{d-1} & \cdots & T^d_{n+1}\end{pmatrix}
$$
so that
$$
H_{n+1}=\E\left[D_{n+1}\,|\,\F_{n}\right]=\E\,D_{n+1}=H=\begin{pmatrix}p^1 & \frac{q^2}{d-1} & \cdots & \frac{q^d}{d-1} \cr \cr \frac{q^1}{d-1} & p^2 & \cdots & \frac{q^d}{d-1} \cr \vdots & \vdots & \ddots & \vdots \cr \frac{q^1}{d-1} & \frac{q^2}{d-1} & \cdots & p^d\end{pmatrix},$$
where $q^i=1-p^i$, $1\leq i\leq d$.  Moreover, $H$ is $\R$-diagonalizable  since its  transpose is  obviously a reversible (stochastic) matrix with respect to its invariant probability measure  $v^*$\footnote{By reversible we mean that $H{\rm diag}(v^*)={\rm diag}(v^*)H^t$. So $H$ is diagonalizable, since it is auto-adjoint with respect to the inner product induced by ${\rm diag}(v^*)$.}, given for this model by 
$$
v^{*i}=\frac{1}{q_i\sum_{j=1}^d1/q^j}, \quad 1\leq i\leq d.
$$
The strong consistency has been first established in~\cite{AthKar2}, then redone in~\cite{BaiHu2}. It follows from Theorem~\ref{Thm1} as well. If $ \lambda_{{\rm max}}<1/2$, the asymptotic normality 
$$
\frac{Y_n-nv^*}{\sqrt{n}}= \sqrt{n} \Big(\frac{Y_n}{n}-v^*\Big)\overset{{\cal L}}{\longrightarrow}{\cal N}(0,\Sigma) \quad \mbox{as $n\to+\infty$}
$$
results from Theorem 3.2 in~\cite{BaiHu2} and from Theorem~\ref{Thm2} of this paper.  Therefore, the other types of rate, depending on $\lambda_{{\rm max}}$, hold (since $r_n\equiv 0$). However, using Theorem~\ref{Thm2} we obtain a joint $CLT$ for $(\widetilde{Y}_n,\widetilde{N}_n)$. 
Note that if $p^i>p^j$, then $v^{*i}>v^{*j}$. Hence the components $v^{*i}$ are ordered according to the increasing efficiency $p^i$ of the treatments.  Furthermore, it is clear that, if $p^i\uparrow 1$ and all other probabilities $p^j$ stand still,  then 
\[
\lim_{p^i\to 1} v^{*j}=\delta_{ij}.
\]
Consequently, since $v^{*i}$ is the asymptotic probability of assigning treatment $i$ to a patient,  the procedure asymptotically allocates more patients to the most efficient treatment(s).  Following the practitioners, the fact that a marginal allocation of less efficient treatments  is preserved  is justified by some comparison matter.

However this model only takes into account in the addition rule matrix $D_n$ the response of the $n^{th}$ patient without considering the ones of past patients. This led the author to introduce \cite{BaiHuShe} a new model based on statistical observations of the efficiency of the assigned treatments to all past patients.


\subsection{The Bai-Hu-Shen {\em GFU} Model}
\label{trois2}

\medskip
 We consider now the model introduced in~\cite{BaiHuShe} (and considered again in~\cite{BaiHu2}) where $(T_n^i)_{n\geq1}$,$1\leq i\leq d$, are $d$ independent sequences of i.i.d. $\{0,1\}$-valued Bernoulli trials satisfying~(\ref{ProbaSucces}) and the filtration $(\F_n)_{n\geq0}$ is defined as in the previous section.
Let $N_n=(N_n^1,\ldots,N_n^d)^t$ and $S_n=(S_n^1,\ldots,S_n^d)^t$, where $N_n^i=N_{n-1}^i+X_n^i$, $n\geq1$, still denotes the number of times the $i^{th}$ treatment is selected among the first $n$ stages and 
$$
S_n^i=S_{n-1}^i+T_n^iX^i_n, \quad n\geq1,
$$
denotes the {\em number of successes} of the $i^{th}$ treatment among these $N_n^i$ trials, $i=1,\ldots,d$. However, to avoid degeneracy of the procedure, we will make the following initialization assumption
\[
N^i_0=1, \; S^i_0=1,\quad i=1,\dots,d
\]
which makes the above interpretation of these quantities correct ``up to one unit".

\paragraph{Remark.} Like with the Wei model, we can simply assume that $T^i_n$ is a $\{0,1\}$-valued efficiency indicator. \\

Define $\Pi_n=(\Pi_n^1,\ldots,\Pi_n^d)^t$,  where $\Pi_n^i=\frac{S_n^i}{N_n^i}$, $i=1,\ldots,d$. In \cite{BaiHuShe} the authors consider the following addition rule matrices,
$$D_{n+1}=\begin{pmatrix}T_{n+1}^1 & \frac{\Pi_n^1(1-T_{n+1}^2)}{\sum_{j\neq2}\Pi_n^j} & \cdots & \frac{\Pi_n^1(1-T_{n+1}^d)}{\sum_{j\neq d}\Pi_n^j} \cr  &  &  &  \cr \frac{\Pi_n^2(1-T_{n+1}^1)}{\sum_{j\neq1}\Pi_n^j} & T_{n+1}^2 & \cdots & \frac{\Pi_n^2(1-T_{n+1}^d)}{\sum_{j\neq d}\Pi_n^j} \cr \vdots & \vdots & \ddots & \vdots \cr \frac{\Pi_n^d(1-T_{n+1}^1)}{\sum_{j\neq1}^d\Pi_n^j} & \frac{\Pi_n^d(1-T_{n+1}^2)}{\sum_{j\neq2}^d\Pi_n^j} & \cdots & T_{n+1}^d \end{pmatrix},$$
$i.e.$ at stage $n+1$, if the response of the $j^{th}$ treatment is a success, then one ball of type $j$ is added in the urn. Otherwise, $\frac{\Pi_n^i}{\sum_{k\neq j}\Pi_n^k}$ (virtual) balls of type $i$, $i\neq j$, are added. This addition rule matrix clearly satisfies {\bf (A1)}-$(i)$ and {\bf (A2)}. Then, one easily checks that the generating matrices are given by
$$
H_{n+1}=\E\left[D_{n+1}\,|\,\F_n\right]=\begin{pmatrix}p^1 & \frac{\Pi_n^1(1-p^2)}{\sum_{j\neq2}\Pi_n^j} & \cdots & \frac{\Pi_n^1(1-p^d)}{\sum_{j\neq d}\Pi_n^j} \cr  &  &  &  \cr \frac{\Pi_n^2(1-p^1)}{\sum_{j\neq1}\Pi_n^j} & p^2 & \cdots & \frac{\Pi_n^2(1-p^d)}{\sum_{j\neq d}\Pi_n^j} \cr \vdots & \vdots & \ddots & \vdots \cr \frac{\Pi_n^d(1-p^1)}{\sum_{j\neq1}\Pi_n^j} & \frac{\Pi_n^d(1-p^2)}{\sum_{j\neq2}\Pi_n^j} & \cdots & p^d \end{pmatrix}
$$
and satisfy {\bf (A1)}-$(ii)$.
As soon as $Y_0\in\R_+^d\setminus\{0\}$, $H_n\overset{a.s.}{\longrightarrow}H$ (see Lemma~\ref{lemme} below or~\cite{BaiHuShe} when $Y_0\in(0,\infty)^d$) where 
$$
H=\begin{pmatrix}p^1 & \frac{p^1(1-p^2)}{\sum_{j\neq2}p^j}  & \cdots & \frac{p^1(1-p^d)}{\sum_{j\neq d}p^j} \cr  &  &  & \cr \frac{p^2(1-p^1)}{\sum_{j\neq1}p^j} & p^2 & \cdots &\frac{p^2(1-p^d)}{\sum_{j\neq d}p^j} \cr \vdots & \vdots & \ddots & \vdots \cr \frac{p^d(1-p^1)}{\sum_{j\neq1}p^j} & \frac{p^d(1-p^2)}{\sum_{j\neq2}p^j} & \cdots & p^d\end{pmatrix}.
$$
The matrix $H$ is clearly irreducible since $0<p^i<1$, $1\leq i\leq d$, so that Assumption {\bf (A3)} is satisfied. The normalized maximal eigenvector $v^*$ (associated to the eigenvalue 1) is given by
\[
v^{*i} = \frac{p_i\,\sum_{k\neq i}p^k}{(1-p_i)\sum_{1\le j\le d}\frac{p^j}{1-p^j}\sum_{k\neq j}p^k},\quad i=1,\ldots,d.
\]
 In the next section, devoted to rates, we will use again the fact that $H$ is $\R$-diagonalizable, still because its transpose is reversible with respect to its invariant distribution $v^*$. Then calling upon Theorem~\ref{Thm1} (or following the direct proof from~\cite{BaiHuShe}) we obtain
\begin{equation}\label{StrongCons}
\widetilde{Y}_n=\frac{Y_n}{n}\overset{a.s.}{\underset{n\rightarrow+\infty}{\longrightarrow}}v^* \quad\mbox{and}\quad \widetilde{N}_n=\frac{N_n}{n}\overset{a.s.}{\underset{n\rightarrow+\infty}{\longrightarrow}}v^*.
\end{equation}
Note that if $p^i>p^j$, $\frac{p^i}{p^j}\frac{\sum_{k\neq i}p^k}{\sum_{k\neq j}p^k}>1$ and $\frac{1-p^j}{1-p^i}>1$ so that $v^{*i}>v^{*j}$. Hence the entries $v^{*i}$ are ordered according to the increasing efficiency $p^i$ of the treatments. This model can be considered as more ethical than the Wei model since a better treatment will be administrated to more patients. Indeed, when $d>2$, for any $i\neq j$, $1\leq i,j\leq d$, if $p^i>p^j$,
$$\frac{v^{*i}_{BHS}}{v^{*j}_{BHS}}>\frac{v^{*i}_{W}}{v^{*j}_{W}}>1$$
(when $d=2$ both matrices $H$ coincide).

\paragraph{Remark.} Note that in that model the ``balls'' in the urn become virtual since there exists no $N\in\N$ such that, for every $n\geq1$, $ND_n\in{\cal M}_d(\N)$. \\

%

\subsection{Asymptotic normality for multi-arm clinical trials for the BHS \texorpdfstring{$GFU$}{GFU} model}
\label{trois3}

In \cite{BaiHuShe} in order to derive a $CLT$, not with the bias $\E Y_n$ but with $nv^*$, from their own general asymptotic normality result (which statement is similar to Theorem~\ref{Thm2}) the authors need to fulfill the following convergence rate assumption for $H_n$
\begin{equation}\label{CondBaiHu}
	\sum_{n\geq1}\frac{\left\|H_n-H\right\|_{\infty}}{\sqrt{n}}<+\infty 
\end{equation}
where $\left\|\cdot\right\|_{\infty}$ is the norm on $L^{\infty}_{\R^{d\times d}}(\P)$. In \cite{BaiHuShe}, an $a.s.$ rate of decay $|\!|\!|H_n-H|\!|\!||_{\infty}=o(n^{-\frac{1}{4}})$ is shown for this model which is clearly not fast enough to fulfill (\ref{CondBaiHu}).

However, by enlarging the dimension of the structure process of the procedure by considering the $3d$-dimensional $\R^d\times{\cal S}\times[0,1]^d$-valued random sequence
$$\widetilde{\theta}_n=\begin{pmatrix}\widetilde{Y}_n\cr\widetilde{N}_n \cr\widetilde{S}_n\end{pmatrix}\quad\mbox{where}\quad\widetilde{S}_n=\frac{S_n}{n},\quad n\geq1,$$
we will establish that a $CLT$ does hold for the BHS $GFU$ model.

The first step is to notice that the generating matrix $H_{n+1}$ can may be written as a function depending on $\widetilde{S}_n$ and $\widetilde{N}_n$, $i.e.$ $H_{n+1}=\Phi(\widetilde{S}_n,\widetilde{N}_n)$, where $\Phi:\R_+^d\times(0,\infty)^d\rightarrow{\cal M}_{d}(\R)$ is a differentiable function defined by
$$
\Phi(s,\nu)=\left(\Phi^{ij}(s,\nu)\right)_{1\leq i,j\leq d} \quad \mbox{where} \quad 
\left\{
\begin{array}{ll}
\Phi^{ii}(s,\nu)=p^i & 1\leq i\leq d\\
\Phi^{ij}(s,\nu)=\frac{s^i/\nu^i}{\sum_{k\neq j}s^k/\nu^k}\,q^j & 1\leq i,j\leq d,\, i\neq j.\\
\end{array}
\right.
$$
Then the following strong consistency and $CLT$ hold for $(\widetilde{\theta}_n)_{n\geq1}$.

\begin{theo}\label{Thm3}
Assume that $Y_0\in\R_+^d\setminus\{0\}$. We still denote by $\lambda_{\max}$ the highest eigenvalue of $H$ apart from $1$.

\smallskip
\noindent$(a)$ If $ \lambda_{\max} = \max\left({\rm Sp}(H)\setminus\{1\}\right)\subset\left(-\infty,\frac12\right)$, then
$$
\widetilde{\theta}_n\overset{a.s.}{\underset{n\rightarrow+\infty}{\longrightarrow}}  \widetilde{\theta}^*\in{\cal S}^2\times[0,1]^d \quad\mbox{ and }\quad \sqrt{n}\left(\widetilde{\theta}_n-\widetilde{\theta}^*\right)\overset{\L}{\underset{n\rightarrow+\infty}{\longrightarrow}}{\cal N}\left(0,\widetilde{\Sigma}\right),
$$
where
$$\widetilde{\theta}^*:=\left(v^*,v^*,{\rm diag}(p)v^*\right)^t,\quad
\widetilde{\Sigma}=\int_0^{+\infty}e^{-u\left(D\widetilde{h}(\widetilde{\theta}^*)-\frac{I_{3d}}{2}\right)}\widetilde \Gamma e^{-u\left(D\widetilde{h}(\widetilde{\theta}^*)-\frac{I_{3d}}{2}\right)^t}du
$$
with 
$$
\widetilde{\Gamma}=\begin{pmatrix}\displaystyle\sum_{k=1}^d v^{*k}C^k-v^*(v^*)^t & H\left({\rm diag}(v^*)-v^*(v^*)^t\right) & \left({\rm diag}(v^*)-v^*(v^*)^t\right){\rm diag}(p) \cr & & & \cr \left({\rm diag}(v^*)-v^*(v^*)^t\right)^tH^t & {\rm diag}(v^*)-v^*(v^*)^t & \left({\rm diag}(v^*)-v^*(v^*)^t\right){\rm diag}(p) \cr & & & \cr {\rm diag}(p)\left({\rm diag}(v^*)-v^*(v^*)^t\right)^t & {\rm diag}(p) \left({\rm diag}(v^*)-v^*(v^*)^t\right)^t& {\rm diag}(p)\left(v^*-v^*v^{*t}{\rm diag}(p)\right)\end{pmatrix}
$$
where $C^k=(C^k_{ij})_{1\leq i,j\leq d}$, $1\leq k\leq d$, are $d\times d$ positive definite matrices with
$$
C^k_{ij}=\frac{p^ip^j(1-p^k)}{\left(\sum_{\ell\neq k}p^\ell\right)^2}\mathds{1}_{\{i,j\neq k\}}+p^k\mathds{1}_{\{i=j=k\}},$$
and $D\widetilde{h}(\widetilde{\theta}^*)$ is an $\R$-diagonalizable matrix reading 
$$
D\widetilde{h}(\widetilde{\theta}^*)=\begin{pmatrix} I_d-H+v^*\mathbf{1}^t & -\frac{\partial}{\partial \nu}\left(\Phi(s,\nu)y\right)_{|\widetilde{\theta}=\widetilde{\theta}^*} & -\frac{\partial}{\partial s}\left(\Phi(s,\nu)y\right)_{|\widetilde{\theta}=\widetilde{\theta}^*} \cr \cr v^*\mathbf{1}^t-I_d & I_d & 0_{{\cal M}_{d}(\R)} \cr \cr \mbox{diag}(p)\left(v^*\mathbf{1}^t-I_d\right) & 0_{{\cal M}_{d}(\R)} &  I_d \end{pmatrix}.
$$
so that $\mbox{Sp}(D\widetilde{h}(\widetilde{\theta}^*)_{\left|{\cal V}_0^2\times \R^d\right.})=\mbox{Sp}((I_d-H)_{\left|{\bf 1}^{\bot}\right.})= \{1-\lambda,\, \lambda\!\in {\rm Sp}(H)\setminus\{1\}\}\subset\R$ (see Appendix~\ref{app:C}).

\medskip
\noindent $(b)$  If $\lambda_{{\rm max}}=1/2$, then, $\widetilde{\theta}_n \to \widetilde{\theta}^*$ a.s. and 
$$
\sqrt{\frac{n}{\log n}}\left(\widetilde{\theta}_n-\widetilde{\theta}^*\right)\overset{\L}{\underset{n\rightarrow+\infty}{\longrightarrow}}{\cal N}\left(0,\widetilde{\Sigma}\right)\mbox{ with }\widetilde{\Sigma}=\lim_{n\to+\infty}\frac{1}{\log n}\int_0^{\log n}e^{-u\left(D\widetilde{h}(\widetilde{\theta}^*)-\frac{I_{3d}}{2}\right)^t}\Gamma e^{-u\left(D\widetilde{h}(\widetilde{\theta}^*)-\frac{I_{3d}}{2}\right)}du.
$$
$(c)$ If $\lambda_{{\rm max}}>1/2$, then $n^{1-\lambda_{{\rm max}}}\left(\widetilde{\theta}_n-\widetilde{\theta}^*\right)$ $a.s.$ converges as $n\to+\infty$ towards a finite random variable.
\end{theo}

%
\noindent {\bf Proof.} {\sc Step~1} ({\em Strong consistency}). We will show with Lemma~\ref{lemme} that $\widetilde{S}_n\overset{a.s.}{\underset{n\rightarrow+\infty}{\longrightarrow}}{\rm diag}(p)v^*$ and we will deduce that $H_n\overset{a.s.}{\underset{n\rightarrow+\infty}{\longrightarrow}}H$, $i.e.$  Assumption~{\bf (A3)} holds. As we have already checked that Assumptions~{\bf (A1)}-$(i)$-$(ii)$ and {\bf (A2)} are satisfied, then by only adding {\bf (A1)}-$(iii)$ we use Theorem~\ref{Thm1} to prove that $\widetilde{\theta}_n\overset{a.s.}{\underset{n\rightarrow+\infty}{\longrightarrow}}  \widetilde{\theta}^*$.

\begin{lem}\label{lemme} If the assumption~(\ref{dynamic}) holds and $Y_0\in\R^d_+\setminus\{0\}$, then, 
\[
 \Pi_n \stackrel{a.s.}{\longrightarrow} p=(p^1,\ldots,p^d)\quad\mbox{ as $n\to +\infty$}
\] 
so that Assumption~(\ref{A3}) holds $i.e.$ $\displaystyle H_n \overset{a.s.}{\underset{n\rightarrow+\infty}{\longrightarrow}} H$. 
\end{lem} 

\paragraph{Remark.}If we assume that $Y_0^i>0$, $1\leq i\leq d$, then we can prove that $\lim_n N^i_n=+\infty $ $a.s.$, $1\leq i\leq d$, faster than below by using that $Y_n^i\geq Y_0^i$, $1\leq i\leq d$, $n\geq1$. The following proof considers the more general case where $Y_0\in\R^d_+\setminus\{0\}$. \\

\noindent {\bf Proof of Lemma~\ref{lemme}.} {\em Step~1.} It follows from the dynamics~(\ref{dynamic}) and the definitions of $D_{n+1}$ and $H_{n+1}$ that, for every $n\ge 0$, $\Tr (Y_n)= \Tr (Y_0)+n$ and that, for every $i\!\in \{1,\ldots,d\}$, 
\[
Y^i_{n+1}=Y^i_{n}+\sum_{j=1}^d H^{ij}_{n+1}  \frac{Y^i_n}{\Tr(Y_n)} +\Delta M^i_{n+1}
\]
where $(\Delta M^i_{n})_{n\ge 1}$ is a sequence of martingale increments satisfying $\sup_n \E \left[|\Delta M^i_n|^2\left.\right|\F_{n-1}\right]<+\infty$ since the addition rule matrices satisfy~(\ref{A2}). Now using that $S^i_0=N^i_0=1$ by convention, one derives that
\[
\forall\, i\neq j,\; H^{ij}_{n+1} \ge  \frac{\kappa_0}{n}, \;\mbox{ with } \;\kappa_0 =\frac{1}{2d} \min_{1\le i\le d}\big(p^i, 1-p^i\big)>0
\] 
so that, using that $H^{ii}_{n+1} = p^i$,  there exists a deterministic integer $n_0$ such that for every $n\ge n_0$, 
\begin{eqnarray*}
Y^i_{n+1}& \ge& \Big (1+\frac{p_i}{n}-\frac{\kappa_0}{\Tr(Y_n)}\Big)Y^i_n+  \frac{\kappa_0}{n} +\Delta M^i_{n+1}\\
&\ge &\Big (1+\frac{p_i}{2\Tr(Y_n)}\Big)Y^i_n+\frac{\kappa_0}{n} +\Delta M^i_{n+1}.
\end{eqnarray*}
Standard computations show that, setting $a^i_n = \prod_{k=n_0}^{n-1} (1+\frac{p_i}{2\Tr(Y_n)}\Big)$, $i=1,\ldots,d$, 
\[
\forall\,n\ge n_0,\quad \frac{Y^i_n}{a ^i_n} \ge \frac{Y^i_{n_0}}{a^i_{n_0}}+\sum_{k=n_0+1}^n  \frac{\kappa_0}{a ^i_{k}} +\sum_{k=n_0+1}^n \frac{\Delta M^i_k}{a^i_k}
\]
Since there exists $\kappa_1$, $\kappa_2>0$ such that $\kappa_1  n^{\frac{p^i}{2}} \le a ^i_n \le \kappa_2 n^{\frac{p^i}{2}}$, one has
\[
\forall\,  \eta>0, \quad \sum_{k=n_0+1}^n \frac{\Delta M^i_k}{a^i_k} = o\big(n^{\frac{1-p^i+\eta}{2}}\big).
\]
Finally, there exists a positive real constant $c'$ such that, for every $i=1,\ldots,d$,  
\[
Y^i_n \ge c' n^{\frac{p^i}{2}}\sum_{k=n_0+1}^n k^{-\frac{p^i}{2} }+ o\big(n^{\frac{1+\eta}{2}}\big)
\] 
so that 
\[
\forall\, i\!\in \{1,\ldots,d\},\quad \liminf_n  \widetilde Y^i_n\, \ge c'  \int_0^1 u^{-\frac{p^i}{2}}du >0
\]
and, as a consequence, $\sum_{n\ge 1} \widetilde Y^i_n =+\infty$ $a.s.$ Now using that  for every $i=1,\ldots,d$, 
$$
N^i_n=\sum_{k=1}^n \mathds{1}_{\{X_k =e ^i\}}\quad\mbox{ and }\quad \P(X_n =e ^i\,|\, \F_{n-1}) = \widetilde Y^i_{n-1}\left(1-\frac{\Tr(Y_0)}{\Tr(Y_{n-1})}\right),\quad n\ge 1,
$$  
we get by the conditional Borel-Cantelli Lemma that $
N^i_{\infty}=\lim_n N^i_n=+\infty $ $a.s.$

\medskip
\noindent {\em Step~2.}  
%
First we note that
\[
\Pi^i_n =\frac{\sum_{k=1}^n T^i_{k}\Delta N^i_k}{N ^i_n}
\]
and we introduce the sequence $(\widetilde \Pi_n)_{n\ge 1}$ defined by 
\[
\widetilde \Pi^i_n = \sum_{k=1}^n (T^i_{k}-p^i)\frac{\Delta N^i_k}{N^i_{k-1}+1},\quad n\ge 1.
\]
It is an ${\cal F}_n$-martingale since, $T^i_k$ being independent of ${\cal F}_{k-1}$ and $X_k$,
\[
\E\Big( (T^i_{k}-p^i)\Delta N^i_k\,|\, {\cal F}_{k-1}\Big)= \E (T^i_k-p^i) \P(X_k= e^i\,|\, {\cal F}_{k-1})=0.
\]
It has bounded increments since $|T^i_{k}-p^i|\le 1$ and
\[
\langle \widetilde \Pi^i\rangle_n  \le  \sum_{k=1}^{n} \frac{\E ((\Delta N^i_k)^2\,|\, {\cal F}_{k-1})}{(N^i_{k-1}+1)^2}. 
\]
It follows, using $(\Delta N^i_k)^2=\Delta N^i_k$, that, for every $n\ge 1$, 
\begin{eqnarray*}
\E \langle \widetilde \Pi^i\rangle_n& \le &\E\Big(\sum_{k=1}^{n} \frac{\Delta N^i_k}{(N^i_{k-1}+1)^2}\Big)\le \E\Big(\sum_{k=1}^{n} \frac{\Delta N^i_k}{N^i_{k-1}N^i_k}\Big)\le \frac {1}{N^i_0}=1.
\end{eqnarray*}
Consequently $\widetilde \Pi^i_n \to \widetilde \Pi^i_{\infty}\!\in L^1(\P)$ $a.s.$ as $n\to +\infty$. This  in turn implies by Kronecker's Lemma that 
$$
\Pi^i_n\stackrel{a.s.}{\longrightarrow} p^i\quad\mbox{ as }\quad n\to +\infty
$$ 
since $N^i_n \to +\infty$ by the first step.~\hfill$\cqfd$

\bigskip It follows from the lemma and Theorem~\ref{Thm1}  that  $(\widetilde Y_n,\widetilde N_n)\to (v^*,v^*)$. Furthermore ${\rm diag}(\widetilde S_n) ={\rm diag}( Q_n) \widetilde N_n\to {\rm diag}(p)v^{*}= u^{*}$  so that $\widetilde \theta_n \to \widetilde \theta^*$ as $n\to+\infty$.

\medskip
\noindent {\sc Step~2} ({\em Asymptotic normality}). We will show now  that $(\widetilde \theta_n)_{n\ge 1}$ satisfies an appropriate recursion to apply Theorem~\ref{ThmCLT}$(a)$ (standard $CLT$). First, we write a recursive procedure for $\widetilde{S}_n$. 
Having in mind that $S_n=1+\sum_{1\le k\le n}\mbox{diag}(T_k)X_k$, we get
\begin{eqnarray}
\widetilde{S}_{n+1}&=&\widetilde{S}_n-\frac{1}{n+1}\left(\widetilde{S}_n-\mbox{diag}(T_{n+1})X_{n+1}\right) \nonumber\\
									 &=&\widetilde{S}_n-\frac{1}{n+1}\left(\widetilde{S}_n-\mbox{diag}(p)\frac{\widetilde{Y}_n}{\Tr(\widetilde{Y}_n)}\right)+\frac{1}{n+1}\Delta \widehat{M}_{n+1} \nonumber\\
									 &=&\widetilde{S}_n-\frac{1}{n+1}\left(\widetilde{S}_n-\mbox{diag}(p)(2-\Tr(\widetilde{Y}_n))\widetilde{Y}_n\right)+\frac{1}{n+1}\left(\Delta \widehat{M}_{n+1}+\widehat{r}_{n+1}\right)\label{AS_Stilde}
\end{eqnarray}
$$\mbox{where  }\quad 
\Delta \widehat{M}_{n+1}:=\mbox{diag}(T_{n+1})X_{n+1}-\E\left[\mbox{diag}(T_{n+1})X_{n+1}\,|\,\F_n\right]=\mbox{diag}(T_{n+1})X_{n+1}-\mbox{diag}(p)\displaystyle\frac{\widetilde{Y}_n}{\Tr(Y_n)}
$$
is an $\F_n$-martingale increment and $\widehat{r}_{n+1}=\mbox{diag}(p)\frac{\left(\Tr(\widetilde{Y}_n)-1\right)^2}{\Tr(\widetilde{Y}_n)}\widetilde{Y}_n$. 
Then we rewrite  the dynamics satisfied by  $\widetilde{Y}_n$ as follows
\begin{equation}\label{AS_Ytilde}
\widetilde{Y}_{n+1}=\widetilde{Y}_n-\frac{1}{n+1}\left(I_d-(2-\Tr(\widetilde{Y}_n))H_{n+1}\right)\widetilde{Y}_n+\frac{1}{n+1}\left(\Delta M_{n+1}+\check{r}_{n+1}\right),
\end{equation}
where $\check{r}_{n+1}:=\displaystyle\frac{\left(\Tr(\widetilde{Y}_n)-1\right)^2}{\Tr(\widetilde{Y}_n)}H_{n+1}\widetilde{Y}_n$. 
%
\noindent Finally, we get the following recursive procedure for $\widetilde{\theta}_n$
$$
\widetilde{\theta}_{n+1}=\widetilde{\theta}_n-\frac{1}{n+1}\widetilde{h}(\widetilde{\theta}_n)+\frac{1}{n+1}\left(\Delta \widetilde{{\bf M}}_{n+1}+\widetilde{R}_{n+1}\right), \quad n\ge1,
$$
where, for every  $ \widetilde{\theta}=(y,\nu, s)^t\!\in\R^{3d}_+$,
$$
 \widetilde{h}(\widetilde{\theta}):=\!\begin{pmatrix} (I_d-(2-\Tr(y))\Phi(s,\nu))y \cr \nu-(2-\Tr(y))y \cr s-(2-\Tr(y))\mbox{diag}(p)y\end{pmatrix}\!,\;
\Delta \widetilde{{\bf M}}_{n+1}:=\begin{pmatrix} \Delta M_{n+1} \cr  \Delta \widetilde{M}_{n+1} \cr \Delta \widehat{M}_{n+1}\end{pmatrix}\,\mbox{ and }\,\widetilde{R}_{n+1}:=\!\begin{pmatrix} \check{r}_{n+1} \cr  \widetilde{r}_{n+1} \cr  \widehat{r}_{n+1}\end{pmatrix}.
$$
Let us check that the addition rule matrices satisfy {\bf (A4)}. For every $j\in\{1,\ldots,d\}$, let set $C_n^j=\E\left[D^{\cdot j}_{n+1}(D^{\cdot j}_{n+1})^t\left|\right.\F_n\right]$. We have that
\begin{eqnarray*}
(C_n^j)_{ii'}&=&\E\left[D^{ij}_{n+1}(D^{i'j}_{n+1})^t\left|\right.\F_n\right]\\
             &=&\frac{Q_n^iQ_n^{i'}}{\left(\sum_{k\neq j}Q_n^k\right)^2}\E\left[(1-T_{n+1}^j)^2\left|\right.\F_n\right]\mathds{1}_{\{i,i'\neq j\}}+\E\left[(T_{n+1}^j)^2\left|\right.\F_n\right]\mathds{1}_{\{i=i'=j\}}
\end{eqnarray*}
because $T_{n+1}^j(1-T_{n+1}^j)=0$. Then owing to Lemma~\ref{lemme}, $C_n^j\overset{a.s.}{\underset{n\to+\infty}{\longrightarrow}}C^j$ with $$C^j_{ii'}=\frac{p^ip^{i'}(1-p^j)}{\left(\sum_{k\neq j}p^k\right)^2}\mathds{1}_{\{i,i'\neq j\}}+p^j\mathds{1}_{\{i=i'=j\}}.$$
We can check that $C^j$ is a positive definite matrice. Consequently {\bf (A4)} holds.

\noindent The function $\Phi$ being differentiable at the equilibrium point $\widetilde{\theta}^*$, we have
\begin{equation}\label{eq:Dhthetattilde}
D\widetilde{h}(\widetilde{\theta}^*)=\begin{pmatrix} I_d-H+v^*\mathbf{1}^t & -\frac{\partial}{\partial \nu}\left(\Phi(s,\nu)y\right)_{|\widetilde{\theta}=\widetilde{\theta}^*} & -\frac{\partial}{\partial s}\left(\Phi(s,\nu)y\right)_{|\widetilde{\theta}=\widetilde{\theta}^*} \cr \cr v^*\mathbf{1}^t-I_d & I_d & 0_{{\cal M}_{d}(\R)} \cr \cr \mbox{diag}(p)\left(v^*\mathbf{1}^t-I_d\right) & 0_{{\cal M}_{d}(\R)} &  I_d \end{pmatrix}.
\end{equation}
Elementary though tedious computations show that 
$\frac{\partial}{\partial \nu}\left(\Phi(s,\nu)y\right)_{|\widetilde{\theta}=\widetilde{\theta}^*}=-\frac{\partial}{\partial s}\left(\Phi(s,\nu)y\right)_{|\widetilde{\theta}=\widetilde{\theta}^*}\mbox{diag}(p)$. It follows (see~Appendix~\ref{app:C}) that $D\widetilde h(\widetilde\theta^*)$ is diagonalizable and that   $\mbox{Sp}(D\widetilde{h}(\widetilde{\theta}^*))=\mbox{Sp}(I_d-H)$. Moreover as $v^*\mathbf{1}^tu=v^*\sum_{i=1}^du^i=0$,  $D\widetilde{h}(\widetilde{\theta}^*)$ leaves stable ${\cal V}_0^2\times \R^d$ and its spectrum on this subspace does not contain $1$, hence is equal to~$\mbox{Sp}((I_d-H))_{\mbox{\bf 1}^\perp}= \{1-\lambda, \, \lambda \!\in {\rm Sp}(H), \, \lambda\neq 1\}$.

As for the reminder term $\widetilde R_{n+1}$ we first note that it is ${\cal F}_n$-measurable and reads 
\[
\widetilde R_{n+1}  = \frac{\left(\Tr(\widetilde{Y}_n)-1\right)^2}{\Tr(\widetilde{Y}_n)} \begin{pmatrix}H_{n+1}\widetilde{Y}_n\cr \widetilde{Y}_n\cr \widetilde{Y}_n.\end{pmatrix} 
\]
As $\frac{\widetilde{Y}_n}{\Tr(\widetilde{Y}_n)}$ lies in the simplex, its $\ell^1$-norm ( $\|(u^1,\ldots,u^d\|_{\ell^1}= |u^1|+\cdots|u^d|$) is $1$ and so is the case of $H_{n+1}\widetilde{Y}_n$.  Finally, following th elines of the end of the proof of Theorem~\ref{Thm2}$(a)$. 
$$
\E\big(\|\widetilde R_{n+1}\|^2\mbox{\bf 1}_{\{\|\widetilde \theta_n- \widetilde\theta^*\|\le \varepsilon\}} \big)\le 3d \E\big( (\Tr(\widetilde{Y}_n)-1)^4\mbox{\bf 1}_{\{\|\widetilde \theta_n- \widetilde\theta^*\|\le \varepsilon\}}\big)\le \frac{C_d}{n^{2+\delta}}
$$
At this stage, the proof follows the lines of that of Theorem~\ref{Thm2}: the computation of the covariance matrix $\widetilde{\Gamma}$ and the treatment of the remainder term uses the same tools as before. 
The three results of convergence rate follow from Theorem~\ref{ThmCLT} in the Appendix (given the above rate obtained for the remainder term). The details are left  to the reader.~\hfill$\cqfd$

\bigskip \noindent {\bf Remark.} The asymptotic variances of $\widetilde{Y}_n$ and $\widetilde{N}_n$ in Theorem \ref{Thm3} are different from those in Theorem \ref{Thm2} because the differential matrices $Dh(\theta^*)$ and $D\widetilde{h}(\widetilde{\theta}^*)$ are not the same.


\begin{cor}\label{Cor1}
Under the assumptions of Theorem~\ref{Thm3}, 
\[
\sqrt{n}\big(H_n-H\big)\overset{\cal L}{\underset{n\rightarrow+\infty}{\longrightarrow}} {\cal N}(0; \Gamma_H)
\]
where  $\Gamma_H$ is a $d^2\times d^2$ matrix given by $\Gamma_H = D\Phi (u^*,v^*)  [\widetilde \Sigma_{i+d,j+d}]_{1\le i,j\le 2d}  D\Phi (u^*,v^*)^t$.
\end{cor}

\noindent {\bf Proof.} This is an easy consequence of the so-called $\Delta$-method since
\[
H_n =\Phi(\widetilde S_n,\widetilde N_n)= \Phi(u^*,v^*) + D\Phi(u^*,v^*).(\widetilde S_n-u^*, \widetilde N_n-v^*) + \|(\widetilde S_n-u^*, \widetilde N_n-v^*) \|\varepsilon(\widetilde S_n, \widetilde N_n) 
\]
with $\lim_{y\to (u^*,v^*)} \varepsilon(y)=0$. Consequently
\[
\sqrt{n}\big( H_n-H\big)=  D\Phi(u^*,v^*).(\sqrt{n}(\widetilde S_n-u^*),\sqrt{n}( \widetilde N_n-v^*))  + \varepsilon_{\P}(n)
\]
where $\varepsilon_{\P}(n)$ goes to $0$ in probability (as the product of a tight sequence and an $a.s.$ convergent sequence).  This concludes the proof.~\hfill$\cqfd$

\bigskip
\noindent{\bf Remark.} This corollary shows {\em a posteriori} that it was hopeless to try applying Theorem~\ref{Thm2} in its standard form to establish asymptotic normality for multi-arm clinical trials since the assumption 
{\bf (A5)} cannot be satisfied. Our global $SA$  approach breaks the vicious circle.

\clearpage

\begin{figure}[ht!]
\centering
\includegraphics[width=14cm]{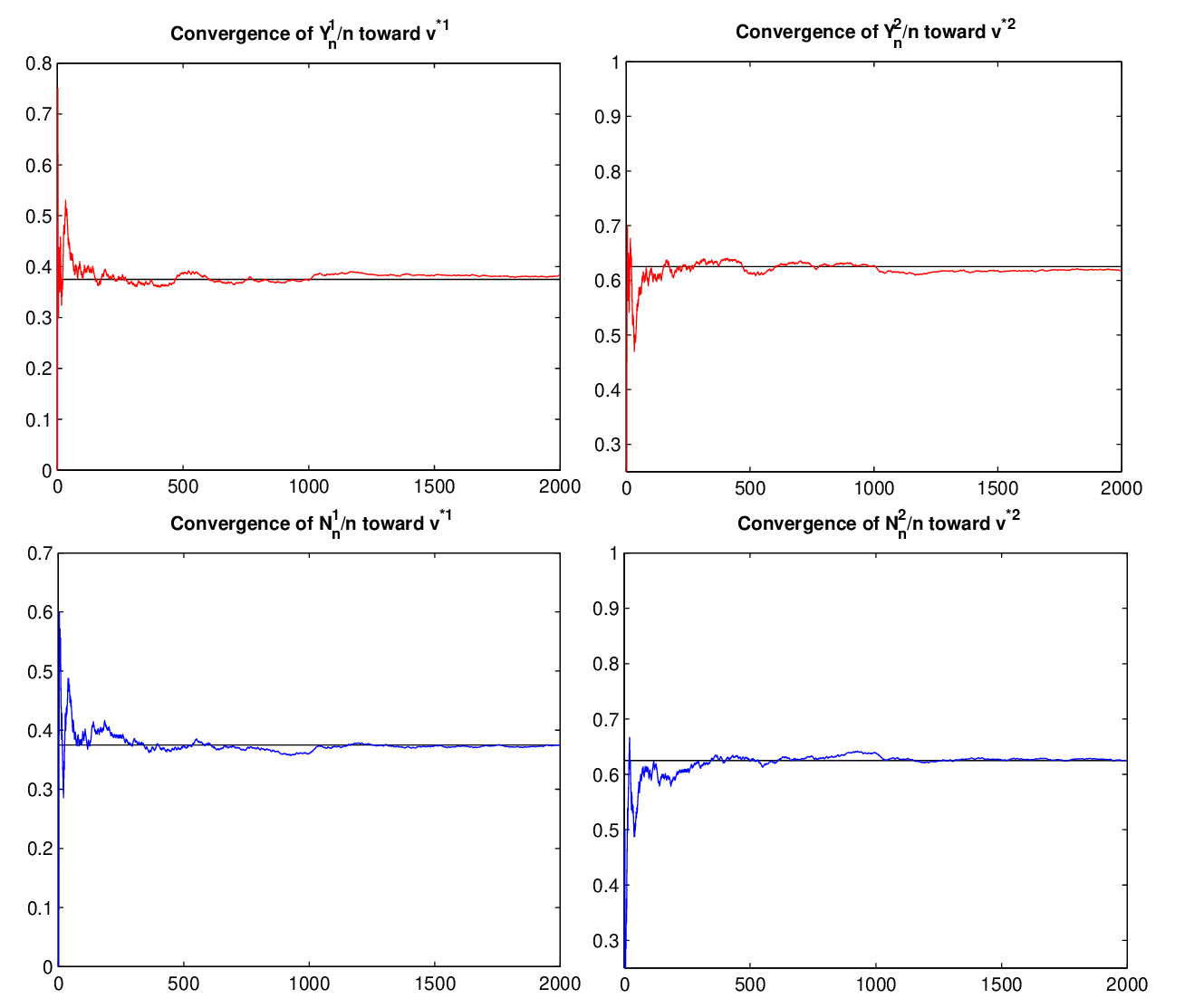}
\caption{Convergence of $\frac{Y_n}{n}$ toward $v^*$ (up-windows) and of $\frac{N_n}{n}$ toward $v^*$ (down-windows): $d=2$, $n=2.10^3$, $p^1=0.5$, $p^2=0.7$, $Y_0=(0.5,0.5)^t$ and $N_0=(1,1)^t$.}
\label{Fig1}
\end{figure}

\bigskip
\noindent {\sc Numerical Example: BHS model}. We consider the case $d=2$, so $v^*$ as the same form as in the example in Subsection~\ref{deux3}. Simulation results are reproduced in Figure~\ref{Fig1}.

\appendix

\begin{center}
\huge{{\bf Appendix}}
\end{center}

\section{Basic tools of Stochastic Approximation}

Consider the following recursive procedure defined on a filtered probability space $(\Omega,{\cal A},(\F_n)_{n\geq0},\P)$
\begin{equation}\label{SAP}
\forall\, n\ge n_0,\quad \theta_{n+1}=\theta_n-\gamma_{n+1}h(\theta_n)+\gamma_{n+1}\left(\Delta M_{n+1}+r_{n+1}\right),
\end{equation}
where $h:\R^d\rightarrow\R^d$ is a locally Lipschitz continuous  function, $\theta_{n_0}$ an $\F_{n_0}$-measurable finite random vector and, for every $n\ge n_0$,   $(\Delta M_{n})$ is a sequence of  $(\F_n)$-martingale increment and $(r_{n})$ is an $(\F_n)$-adapted sequence of remainder terms. 

\begin{theo}[$A.s.$ convergence with $ODE$ method, see $e.g.$~\cite{BMP, Duf2, KusYin, ForPag, Ben}] \label{ThmODE} Assume that $h$ is locally Lipschitz, that
$$
r_n\overset{a.s.}{\underset{n\rightarrow+\infty}{\longrightarrow}}0 \quad \mbox{and} \quad \sup_{n\geq n_0}\E\left[\left\|\Delta M_{n+1}\right\|^2\,|\,\F_n\right]<+\infty \quad a.s.,
$$
and that $(\gamma_n)_{n\geq1}$ is a positive sequence satisfying 
$$
\sum_{n\geq1}\gamma_n=+\infty \quad \mbox{and} \quad \sum_{n\geq1}\gamma_n^2<+\infty.
$$
Then the set $\Theta^{\infty}$ of its limiting values as $n\rightarrow+\infty$ is $a.s.$ a compact connected set, stable by the flow of
$$
ODE_h\equiv\dot{\theta}=-h(\theta).
$$
Furthermore if $\theta^*\in\Theta^{\infty}$ is a uniformly stable equilibrium on $\Theta^{\infty}$ of $ODE_h$, then
$$
\theta_n\overset{a.s.}{\longrightarrow}\theta^*\quad\mbox{ as }\quad  n\to +\infty.
$$
\end{theo}

\paragraph{{\sc Comments.}} By uniformly stable we mean that 
$$
\sup_{\theta\in\Theta^{\infty}}\left|\theta(\theta_0,t)-\theta^*\right| \longrightarrow 0\quad\mbox{as} \quad t\rightarrow+\infty
$$
where $\theta(\theta_0,t)_{\theta_0\in \Theta^{\infty},\,t\in \R_+}$ is the flow of $ODE_h$ on $\Theta^{\infty}$. 

\medskip
We introduce the  $\eta${\em -differentiability} of  the vector field $h$ at $\theta^*$: 
\begin{equation}\label{eta-diff}
h(\theta)=h(\theta^*)+Dh(\theta^*)(\theta-\theta^*)+o\big(\left\|\theta-\theta^*\right\|^{1+\eta}\big) \quad\mbox{as}\quad \theta\to\theta^* \quad\mbox{for some $\eta>0$}.
\end{equation}

\begin{theo}[Rate of convergence see \cite{Duf2} Theorem 3.III.14 p.131~(for $CLT$ see also $e.g.$~\cite{BMP,KusYin})]\label{ThmCLT} Let $\theta^*$ be an equilibrium point of $\{h = 0\}$. Assume that the function $h$ is differentiable at $\theta^*$ and all the eigenvalues of $Dh(\theta^*)$ have positive real parts. Assume that for some $\delta>0$,
\begin{equation}\label{HypDM}
	\sup_{n\geq n_0}\E\left[\left\|\Delta M_{n+1}\right\|^{2+\delta}\,|\,\F_n\right]<+\infty \, a.s., \quad \E\left[\Delta M_{n+1}\Delta M_{n+1}^t\,|\,\F_n\right]\overset{a.s.}{\underset{n\rightarrow+\infty}{\longrightarrow}}\Gamma,
\end{equation}
where $\Gamma$ is a deterministic symmetric definite positive matrix and for an $\epsilon>0$, 
\begin{equation}\label{HypReste}
(n+1)v_n\E\left[\left\|r_{n+1}\right\|^2\mathds{1}_{\{\left\|\theta_n-\theta^*\right\|\leq\epsilon\}}\right]\underset{n\rightarrow+\infty}{\longrightarrow}0,
\end{equation}
where $(v_n)_{n\geq1}$ is a positive sequence. Specify the gain parameter sequence as follows
\begin{equation}\label{HypoPas}
\forall n\geq1, \quad \gamma_n=\frac{1}{n}.
\end{equation}
$(a)$ If $\Re e(\lambda_{{\rm min}})>\frac{1}{2}$, where $\lambda_{{\rm min}}$ denotes the eigenvalue of $Dh(\theta^*)$ with the lowest real part and~\eqref{HypReste} holds with  $v_n=1$, $n\geq1$,  then, the above $a.s.$ convergence is ruled on the convergence set $\{\theta_n\rightarrow\theta^*\}$ by the following Central Limit Theorem
$$
\sqrt{n}\left(\theta_n-\theta^*\right)\overset{{\cal L}}{\underset{n\rightarrow+\infty}{\longrightarrow}}{\cal N}\left(0,\Sigma\right) \quad
\mbox{with} \quad \Sigma:=\displaystyle\int_0^{+\infty} e^{-\left(Dh(\theta^*)^t-\frac{I_d}{2}\right)u}\Gamma e^{-\left(Dh(\theta^*)-\frac{I_d}{2}\right)u}du.
$$

\noindent $(b)$ If $\Re e(\lambda_{\min})=\frac{1}{2}$, $h$ is $\eta$-differentiable at $\theta^*$ with diagonalizable $Dh(\theta^*)$ and~\eqref{HypReste} holds with $v_n=\log n$, $n\geq2$, then
$$
\sqrt{\frac{n}{\log n}}\left(\theta_n-\theta^*\right)\overset{{\cal L}}{\underset{n\rightarrow+\infty}{\longrightarrow}}{\cal N}(0,\Sigma)\mbox{ with }\Sigma= \lim_{T\to+\infty}\frac{1}{T}\int_0^{T}e^{-\left(Dh(\theta^*)^t-\frac{I_d}{2}\right)u} \Gamma e^{-\left(Dh(\theta^*)-\frac{I_d}{2}\right)u}du.
$$
$(c)$  If $\lambda_{{\rm min}}\!\in(0,\frac{1}{2})$, $Dh(\theta^*)$ is as above  and~\eqref{HypReste} holds with $v_n=n^{2 \lambda_{\min}-1+\varepsilon}$, $n\geq1$, for some $\varepsilon>0$, then $n^{\lambda_{\min}}\left(\theta_n-\theta^*\right)$ $a.s.$ converges as $n\to+\infty$ towards a finite random variable.
\end{theo}

\paragraph{Remark.} After this paper was published in Annals of Applied Probability, L.-X. Zhang pointed out  in~\cite{Zhang} a less stringent assumption on $H$ to get  $(b)$ and $(c)$, namely that  all the Jordan blocks of $\lambda_{{\rm min}}$ have order $1$  and, in~$(c)$, that $\lambda_{\min}$ can be replaced {\em mutatis mutandis} by $\Re e(\lambda_{\min})$. When   these orders are not equal to~$1$  (or even in situations when $H$ itself is random), new    rates  are obtained  (see Theorem~2.1 in~\cite{Zhang}). Thus, in item~$(b)$, if $\nu$ denotes the maximum size of Jordan blocks of $\lambda_{\min}$ then $\sqrt{\frac{n}{\log n}}$ should be replaced by $\frac{\sqrt{n}}{(\log n)^{\nu-\frac 12}}$ (and the definition of $\Sigma$ should be modified accordingly  by replacing $1/T$ by $1/T^{2\nu-1}$ in the r.h.s. of its definition).

\medskip  Note that in our examples of applications the matrices $H$ are diagonalizable since $H^t$  always  turns out to be reversible w.r.t. to their invariant distribution.  

\section{On the eigenvalues of the limit generating matrix in the Bai-Hu-Shen model}\label{app:B}

We have seen in Section~\ref{trois3} that  the limit generating matrix $H$ of the $BHS$ model reads
$$
H=\left(p^i\delta_{ij}+\frac{p^i(1-p^j)}{\pi-p^j}(1-\delta_{ij})\right)_{1\leq i,j\leq d} \quad\mbox{where}\quad\pi=\sum_{i=1}^dp^i
$$
and is always diagonalizable since its transpose is reversible with respect to its ``first" eigenvector $v^*$. We propose below another proof when the $p^i$ are pairwise distinct which provides bounds for the eigenvalues. Hence we can give a sufficient condition for having a standard $CLT$ for the urn dynamics. 
\begin{theo}\label{ThmBHS}
The characteristic polynomial of the above $BHS$ generating matrix $H$ is given by $$
\mbox{det}(H-\lambda I_d)=\prod_{i=1}^d\left(p^i(1-a^i)-\lambda\right)+\sum_{i=1}^dp^ia^i\prod_{i\neq j}\left(p^j(1-a^j)-\lambda\right),
$$
where $a^i=\frac{1-p^i}{\pi-p^i}$, $i\in\{1,\ldots,d\}$. In particular, if for every $i\neq j$, $p^i\neq p^j$, then $H$ has pairwise distinct eigenvalues hence it  is diagonalizable. Furthermore the second highest eigenvalue $\lambda_{\max} $ of $H$  satisfies
\[
\lambda_{\max} < \max_{1\le i\le d}\frac{p^i(1-p^i)}{\pi -p^i}.
\]
\end{theo}

\paragraph{Proof.} Setting $D_d(\lambda,p^{1:d},a^{1:d})=\mbox{det}\left[\left(1-\frac{\lambda}{p^i}\right)\delta_{ij}+a^j(1-\delta_{ij})\right]$ implies that $\mbox{det}(H-\lambda I_d)=\prod_{i=1}^dp^iD_d(\lambda,p^{1:d},a^{1:d})$. Moreover, by subtracting the second line to the first one and by developing with respect to the first line, we obtain that
$$D_d(\lambda,p^{1:d},a^{1:d})=\left(1-\frac{\lambda}{p^1}-a^1\right)D_{d-1}(\lambda,p^{2:d},a^{2:d})+a^1\prod_{i=2}^d\left(1-\frac{\lambda}{p^i}-a^i\right).$$
By iteration, we get
$$D_d(\lambda,p^{1:d},a^{1:d})=\prod_{i=1}^d\left(1-\frac{\lambda}{p^i}-a^i\right)+\sum_{i=1}^da^i\prod_{i\neq j}\left(1-\frac{\lambda}{p^j}-a^j\right).$$
Therefore
$$\mbox{det}(H-\lambda I_d)=\prod_{i=1}^d\left(p^i(1-a^i)-\lambda\right)+\sum_{i=1}^dp^ia^i\prod_{i\neq j}\left(p^j(1-a^j)-\lambda\right).$$

\noindent $\rhd$ If for every $i\neq j$, $p^i\neq p^j$ and $\pi\neq1$, then for every $i\neq j$, $p^i(1-a^i)\neq p^j(1-a^j)$. Consequently, there exists a permutation $\sigma\in\Sigma_d$ such that $i\mapsto p^{\sigma(i)}(1-a^{\sigma(i)})$ is increasing. Thus, one checks by considering the function $\lambda\mapsto\frac{\mbox{det}(H-\lambda I_d)}{\prod_{i=1}^d(p^i(1-a^i)-\lambda)}$ that there are $d$ distinct roots for $\mbox{det}(H-\lambda I_d)$ such that $\lambda_i\in(p^{\sigma(i)}(1-a^{\sigma(i)}),p^{\sigma(i+1)}(1-a^{\sigma(i+1)}))$, $i\in\{1,\ldots,d\}$ (with the convention that $p^{\sigma(d+1)}(1-a^{\sigma(d+1)})=+\infty$). Consequently, $H$ has $d$ real distinct eigenvalues. \\

\noindent $\rhd$ If for every $i\neq j$, $p^i\neq p^j$ and $\pi=1$, then 1 is an eigenvalue of $H$ of multiplicity one and 0 of multiplicity $d-1$. It is easy to check that the eigensubspace associated to 0 is of dimension $d-1$. \\

Therefore, if for every $i\neq j$, $p^i\neq p^j$, $H$ is diagonalizable. $\cqfd$

\section{Additional  results}\label{app:C}

In this section we briefly prove that the matrices $Dh(\theta^*)$ are  diagonalizable in both investigated models.
 
\paragraph{Spectrum of $Dh(\theta^*)_{|{\cal V}^2_0}$ in Theorem~\ref{Thm2}.} We aim at proving that, if $H$ is diagonalizable, so is the case of $Dh(\theta^*)_{|{\cal V}^2_0}$.

 We know that $H$ leaves stable $\R v^*$ and ${\cal V}_0$ and $\R^d = \R v^*\oplus {\cal V}_0$. So let $\lambda\!\in {\rm Sp} (H)\setminus\{1\}$ and   $y\!\in E^H_{\lambda}$  (eigenspace of $\lambda$).  Noting  that $E^H_{\lambda}\subset {\cal V}_0$ and that $v^*\mbox{\bf 1}^ty= (\sum_iy_i)v^*=0$, one derives that $\begin{pmatrix}y\cr \frac{y}{1-\lambda}\end{pmatrix}$ is an eigenvector of $Dh(\theta^*)_{|{\cal V}^2_0}$. If $H$ admits a base of eigenvectors $(v^*, y_2, \dots,y_d)$ on $\R^d$, it is clear that if $(\nu_1,\ldots,\nu_{d-1})$ is  basis of ${\cal V}_0$ then 
 \[
 \begin{pmatrix} y_i \cr \cr \frac{y_i}{1-\lambda_i}\end{pmatrix}, \; i=2,\ldots,d,\; \nu_2,\ldots,\nu_{d}, 
 \]
makes up clearly an eigenbasis of ${\cal V}^2_0$  for $Dh(\theta^*)_{|{\cal V}^2_0}$.

\paragraph{Spectrum of $D\widetilde h(\widetilde \theta^*)_{|{\cal V}^2_0}$ in Theorem~\ref{Thm3}.} 
We are interested in the spectrum of $D\widetilde{h}(\widetilde{\theta}^*)_{|{\cal V}_0^2\times \R^d}$ (this vector subspace is left stable by $D\widetilde{h}(\widetilde{\theta}^*)$). Still owing to $v^*\mbox{\bf 1}^ty=0$ for $y\!\in {\cal V}_0$, we derive from~\eqref{eq:Dhthetattilde} that
$$
D\widetilde{h}(\widetilde{\theta}^*)_{|{\cal V}_0^2\times \R^d}=
\begin{pmatrix} I_{\mbox{\bf 1}^{\perp}} -H_{\mbox{\bf 1}^{\perp}} &- B\,{\rm diag}(p) & B\cr \cr 
 I_{\mbox{\bf 1}^{\perp}}  & I_{\mbox{\bf 1}^{\perp}}  & 0_{{\cal M}_{d}(\R)} 
\cr \cr -\mbox{diag}(p) I_{\mbox{\bf 1}^{\perp}}  & 0_{{\cal M}_{d}(\R)} &  I_d \end{pmatrix}.
$$
with $B=  -\frac{\partial}{\partial s}\left(\Phi(s,\nu)y\right)_{|\widetilde{\theta}=\widetilde{\theta}^*} $. Let $y\!\in {\cal V}_0$ be an eigenvector of $H$ with eigenvalue $\lambda \neq 1$. Then, elementary computations show that $(y, \frac{y}{1-\lambda}, \frac{{\rm diag}(p)y}{1-\lambda})^t$ is an eigenvector in ${\cal V}^2_0:times \R^d$ for $D\widetilde{h}(\widetilde{\theta}^*)_{|{\cal V}_0^2\times \R^d}$. As a consequence, if $(y_2,\ldots,y_d)$ is a eigenbasis of $H_{|\mbox{\bf 1}^{\perp}}$ and $(\nu_, \ldots,\nu_{d-1}, e_1,\ldots, e_d)$ denotes a basis of $\{0_{{\cal V}_0}\}\times{\cal V}_0\times\R^d$,  then
\[
\begin{pmatrix}y_i \cr\cr \frac{y_i}{1-\lambda_i}\cr\cr \frac{{\rm diag}(p) y_i}{1-\lambda_i}\end{pmatrix} , \, i=2,\ldots,d,\, \nu_2, \ldots,\nu_{d}, e_1,\ldots, e_d
\]
makes up an eigenbasis of  $D\widetilde{h}(\widetilde{\theta}^*)_{|{\cal V}_0^2\times \R^d}$.
\small

\bibliography{biblio}

%
%
%

\end{document}